\documentclass[11pt,french]{article}
\usepackage{graphics}
\usepackage[french]{babel}
\usepackage[latin1]{inputenc}
\usepackage{amsfonts}

\newtheorem{machin}{Machin}{\bf}{\it}
\newtheorem{theoreme}[machin]{Théorème}{\bf}{\it}
\newtheorem{corollaire}[machin]{Corollaire}{\bf}{\it}
\newtheorem{lemme}[machin]{Lemme}{\bf}{\it}
\newtheorem{propositio}[machin]{Proposition}{\bf}{\it}
{\bf}{\it}
\newtheorem{definitio}[machin]{Définition}{\bf}{\it}

\newcommand{\cad}{c'est-à-dire}
\newcommand{\grad}{\ensuremath{\mathbf{\nabla}}}
\newcommand{\sii}{si, et seulement si}
\newcommand{\tq}{\ensuremath{\;|\;}}
\newcommand{\qed}{\ensuremath{\;\diamondsuit\;}}

\newcommand{\N}{\ensuremath{\mathbb{N}}}
\newcommand{\Z}{\ensuremath{\mathbb{Z}}}

\newcommand{\Zsd}{\ensuremath{\mathbb{Z}_{2}}}

\newcommand{\R}{\ensuremath{\mathbb{R}}}
\newcommand{\Rd}{\ensuremath{\mathbb{R}^{2}}}
\newcommand{\Rt}{\ensuremath{\mathbb{R}^{3}}}

\newcommand{\Rn}{\ensuremath{\mathbb{R}^{n}}}
\newcommand{\C}{\ensuremath{\mathbb{C}}}
\newcommand{\Cd}{\ensuremath{\mathbb{C}^{2}}}

\newcommand{\Cn}{\ensuremath{\mathbb{C}^{n}}}
\newcommand{\CPu}{\ensuremath{\mathbb{CP}^{1}}}
\newcommand{\CPd}{\ensuremath{\mathbb{CP}^{2}}}
\newcommand{\CPdual}{\mbox{\boldmath $H$}}

\newcommand{\CPn}{\ensuremath{\mathbb{CP}^{n}}}

\newcommand{\CPM}{\ensuremath{\mathbb{CP}^{M}}}

\newcommand{\eit}{\ensuremath{e^{i\theta}}}

\newcommand{\espV}{\ensuremath{\mathbf{V}}}

\newcommand{\vecx}{\ensuremath{\mathbf{x}}}
\newcommand{\vecy}{\ensuremath{\mathbf{y}}}

\newcommand{\vecv}{\ensuremath{\mathbf{v}}}
\newcommand{\vecw}{\ensuremath{\mathbf{w}}}

\newcommand{\spbd}{\mbox{\boldmath $s$}}

\newcommand{\grass}{\mbox{\boldmath $G$}}
\newcommand{\Splus}{\ensuremath{S_{+}}}
\newcommand{\Smoins}{\ensuremath{S_{-}}}

\newcommand{\Cinf}{\ensuremath{C^{\infty}}}
\newcommand{\holostar}{\ensuremath{\mathcal{O}^{\ast}}}

\newcommand{\vU}{\ensuremath{U}}
\newcommand{\vV}{\ensuremath{V}}
\newcommand{\vW}{\ensuremath{W}}
\newcommand{\vX}{\ensuremath{X}}

\newcommand{\vY}{\ensuremath{Y}}

\newcommand{\fcX}{\ensuremath{\mathcal{K}_{X}}}
\newcommand{\fcS}{\ensuremath{\mathcal{K}_{\Sigma}}}
\newcommand{\vXz}{\ensuremath{\vX_{0}}}

\newcommand{\vXt}{\ensuremath{\tilde{\vX}}}

\newcommand{\Sigmat}{\ensuremath{\tilde{\Sigma}}}
\newcommand{\Sigmaz}{\ensuremath{\Sigma_{0}}}

\newcommand{\g}{\mbox{$\mathsf{g}$}}
\newcommand{\gp}{\mbox{$\mathsf
g$\textsf{'}}}
\newcommand{\gz}{\mbox{$\mathsf{g}_{0}$}}

\newcommand{\pts}{\mbox{\sl s}}
\newcommand{\ptt}{\mbox{\sl t}}

\newcommand{\ptx}{\mbox{\sl x}}
\newcommand{\pty}{\mbox{\sl y}}
\newcommand{\ptz}{\mbox{\sl z}}

\newcommand{\omegaz}{\ensuremath{\omega_{0}}}
\newcommand{\lambdaz}{\ensuremath{\lambda_{0}}}
\newcommand{\omegap}{\ensuremath{\omega'}}
\newcommand{\Jt}{\ensuremath{\tilde{J}}}
\newcommand{\Jz}{\ensuremath{J_{0}}}

\newcommand{\vectx}{\mbox{\boldmath $x$}}
\newcommand{\vecty}{\mbox{\boldmath $y$}}
\newcommand{\vectn}{\mbox{\boldmath $n$}}
\newcommand{\vectxi}{\mbox{\boldmath $\xi$}}

\newcommand{\vectu}{\mbox{\boldmath $u$}}

\newcommand{\Sigmamth}{\ensuremath{\Sigma_{m}^{(\theta)}}}

\newcommand{\Sigmakth}{\ensuremath{\Sigma_{k}^{(\theta)}}}
\newcommand{\skth}{\ensuremath{s_{k}^{(\theta)}}}
\newcommand{\Akth}{\ensuremath{A_{k}^{(\theta)}}}
\newcommand{\Transv}{\mbox{\boldmath $T$}}
\newcommand{\Hs}{\mbox{{\boldmath $H$}$_{\ens}$}}
\newcommand{\Vtransv}{\mbox{\boldmath $V$}}

\newcommand{\prodint}{\mbox{$\rule{1.5mm}{.15mm}\rule{.15mm}{2mm}\:$}}

\newcommand{\HdZ}[1]{\ensuremath{H_{2}(#1;\Z)}}
\newcommand{\HdR}[1]{\ensuremath{H_{2}(#1;\R)}}

\newcommand{\HDZ}[1]{\ensuremath{H^{2}(#1;\Z)}}

\newcommand{\ens}{\mbox{\scriptsize\sl s}}
\newcommand{\enx}{\mbox{\scriptsize\sl x}}
\newcommand{\eny}{\mbox{\scriptsize\sl y}}
\newcommand{\enz}{\mbox{\scriptsize\sl z}}
\newcommand{\degmetr}{\mbox{\scriptsize\sf g}}
\newcommand{\dey}{\mbox{\scriptsize\bf y}}
\newcommand{\surV}{\mbox{$\scriptstyle{V}$}}
\newcommand{\surX}{\mbox{$\scriptstyle{X}$}}

\begin{document}
\title{Représentants lagrangiens de l'homologie des surfaces projectives complexes}
\author{Daniel Bennequin\\ Université Paris 7 / Denis Diderot\\ Institut de
Mathématiques - UMR CNRS 7586 \and Thanh-Tâm Lê\\ Ecole
Polytechnique}
\date{16 mars 2009}
\maketitle
\begin{abstract}
Using results by Donaldson and Auroux on pseudo-holo\-morphic curves
as well as Duval's rational convexity construction, the paper
investigates the existence of smooth Lagrangian surfaces
representing 2-dimensional homology classes in complex projective
surfaces. We prove that if the projective surface $\vX$ is minimal,
of general type, with uneven geometric genus, and has an effective,
smooth and connected canonical divisor $\mathcal{K}$, then there
exists a non-empty convex open cone in the real 2-dimensional
homology group $\HdR{\vX}$ such that a multiple of every integral
homology class in this cone can be represented by an embedded
Lagrangian surface in $\vX\setminus\mathcal{K}$. A corollary asserts
that such a surface $\vX$ is of simple type in the sense of
Kronheimer and Mrowka.
\end{abstract}

\section{Introduction}
L'ambition de décrire la topologie des surfaces projectives
complexes est ancienne, comme le montrent les travaux fondateurs de
Picard et Lefschetz \cite{Lef24,PiSi}. Elle a connu récemment des
dévelop\-pements spectaculaires, notamment à travers les travaux de
Donaldson, Seiberg-Witten et Taubes. Mais il manque encore, pour les
surfaces complexes, une représentation géométrique aussi claire que
celle qu'on a pour les surfaces de Riemann ou maintenant pour les
variétés de dimension 3.

Les plus gros progrès sur la question sont dus à l'analyse globale,
d'abord développée dans le cadre linéaire par Hodge, Atiyah,
Hirzebruch et Kodaira. Ensuite apparaît l'analyse non-linéaire, en
particulier l'étude des instantons (Atiyah, Hitchin, Donaldson,
Kronheimer et Mrowka) et des monopôles (Seiberg et Witten, Taubes).
Aujourd'hui interviennent les théories des super-cordes et des
D-branes.

Notons que les surfaces projectives complexes sont naturellement des
variétés symplectiques réelles compactes de dimension 4~: on peut
donc parler de leurs sous-variétés symplectiques ou lagrangiennes.
La plupart des avancées géométriques récentes font explicitement
intervenir la géométrie symplectique. D'ailleurs, elles réussissent
à étendre aux variétés symplectiques une partie des résultats connus
dans le cas des surfaces projectives complexes.

Deux techniques fondamentales ont été introduites dans ce sens~:
celle des courbes pseudo-holomorphes de Gromov \cite{Gro} et celle
des sections asymptotiquement holomorphes de fibrés développée par
Donaldson \cite{Don96}. Ces méthodes ont permis d'établir
l'existence de surfaces symplectiques plongées dans les variétés
symplectiques de dimension 4
\cite{Gro,TauGrSW,TauSWGr,Aur97,Aur00,Aur01,DoSm}.

Du côté des surfaces lagrangiennes plongées, on dispose de très peu
de résultats d'existence~; en revanche, les propriétés de rigidité,
d'unicité, ou non, à isotopie près ont été mieux explorées
\cite{Eli90a,Pol,Mok,Sei97,Bir01,BiCi,Auc,Vid}.

Pour étudier le problème de la disjonction des sous-variétés lagrangiennes
(dont le problème d'existence des orbites périodiques de flots hamiltoniens est
un cas particulier), A.~Floer a introduit sa théorie d'homologie
\cite{Flo88,Flo89}. Partant de là, Fukaya \cite{Fuk93} puis Fukaya, Oh, Ohta et
Ono (\textit{cf.} \cite{Fuk03,OhFu}) ont défini une théorie cohomologique
($A_{\infty}$-catégorie), qui a permis à Kontsevich et Soibelman de formuler un
ensemble de conjectures appelées ``conjectures miroir homologiques''
\cite{KoSo01,KoSo06}.

Par exemple, pour les surfaces K3 et, plus généralement, pour les variétés de
Calabi-Yau, ces conjectures proposent une équivalence de catégories entre la
structure $A_{\infty}$ de Fukaya d'une variété et la catégorie dérivée des
faisceaux cohérents d'une autre variété, sa vari\-été ``duale''~:
heuristiquement, certaines sous-variétés lagrangiennes d'une variété
quasi-projective correspondent à des diviseurs holomorphes d'une autre variété
(\textit{cf.} \cite{GrWi}). La symétrie miroir provient des théories quantiques
de champs (super-)conformes issues des théories de super-cordes. Cela reste un
problème de comprendre son extension la plus générale. Auroux, Katzarkov et
Orlov \cite{AKO} étendent la dualité de surfaces à certaines variétés
rationnelles et à certaines fibrations elliptiques.

En un sens, les constructions de notre article suggèrent une
extension de la dualité miroir aux surfaces de type général~:
partant d'une surface projective $\vX$ de type général avec un
diviseur canonique lisse, elles construisent une variété
symplectique $\tilde{X}$ dont les courbes pseudo-holomorphes donnent
des surfaces lagrangiennes de $\vX$. Signalons que, dans sa thèse
\cite{Bah}, Alireza Bahraini procède à une construction analogue
mais où $\tilde{X}$ est munie d'une structure analytique complexe
intégrable singulière le long d'un diviseur.

Dans cet article, nous établissons le premier théorème général
d'existence de surfaces lagrangiennes de genre plus grand que 2.

Nous nous plaçons dans le cas où la surface algébrique complexe
$\vX$ possède une structure kählérienne  $\omegaz$ colinéaire à sa
classe canonique $K$~: c'est ce qui arrive pour les plongements
pluricanoniques des surfaces de type général.

Suivant la démarche initiale de Picard, une fois qu'on dispose d'un
diviseur $\Sigma$ représentant $\omegaz$ (ou un multiple), on peut
considérer le complémentaire de $\Sigma$ dans $\vX$~: c'est une
variété de Stein dont l'homologie a été appelée par Picard
l'``homologie à distance finie'' de $\vX$. Il est très tentant de
chercher à rétracter $\vX\setminus\Sigma$ sur un 2-complexe
cellulaire isotrope. C'est ce que devrait faire le gradient d'une
fonction pluri-sous-harmonique \cite{AnFr} mais, en général, il se
pose un problème de désingularisation. En utilisant la construction
asymptotique de Donaldson, Paul Biran a démontré l'existence d'un
tel complexe avec des intersections transverses. Mais il demeure le
problème de rendre ce complexe le plus explicite et économique
possible.

Par ailleurs, dans l'esprit indiqué par Lefschetz, on aimerait
repré\-senter les classes d'homologie orthogonales à $\omegaz$ (et à
$K$) par des surfaces lagrangiennes dans $\vX\setminus\Sigma$ et
leur trouver les représentants lisses les plus simples possibles (on
sait qu'il arrive que des surfaces lagrangiennes non isotopes soient
homologues (\textit{cf.} \cite{FiSt04,Vid}). Paul Seidel, dans sa
thèse \cite{Sei97}, a construit des représentants lagrangiens
sphériques pour les cycles évanescents~; à l'aide de ces cycles, il
a étudié les groupes de difféomorphismes symplectiques à isotopie
près, puis il a formé une théorie cohomologique ``évanescente'' dans
l'esprit de Fukaya \cite{Sei01a,Sei01b}.

Signalons que le problème de l'existence et de l'unicité d'une
sous-variété lagrangienne avec des conditions homologiques
prescrites a donné lieu à plusieurs études récentes d'un grand
intérêt \cite{ADK,Bir06,BiCi,Alb,Chi,Sei00,Wol}.


\`{A} l'origine du présent travail se trouve un exposé de Daniel
Bennequin (lors du colloque Thom de 1988) qui portait sur
l'existence de squelettes lagrangiens et de représentants
lagrangiens dans les classes d'homologie. Le théorème principal du
présent article faisait partie de la thèse de Lê Thanh-Tâm soutenue
en 2002 à l'Université Paris 7 \cite{Le02}.

Notre travail utilise les résultats de Donaldson et d'Auroux sur les
courbes pseudo-holomorphes \cite{Don96,Don99,Aur97,Aur00,Aur01},
ainsi qu'une construction de Duval sur la convexité rationnelle
\cite{Duv}.

Nous démontrons que si la surface projective $\vX$ est minimale, de
type général, de genre géométrique impair et admet un diviseur
canonique $K$ effectif, lisse et connexe, alors il existe un cône
ouvert non-vide dans l'homologie réelle $H_{2}$ de $\vX\setminus K$
tel que toute classe d'homologie entière dans ce cône admet un
multiple qui se représente par une surface lagrangienne plongée dans
$\vX\setminus K$.

Un corollaire affirme qu'une telle surface $\vX$ est de type simple
au sens de Kronheimer et Mrowka \cite{KrMr95}.

\vspace{5mm}

\noindent \begin{small} \textit{Remerciements.} Les auteurs tiennent
à remercier, pour leur aide et leur participation, Denis Auroux,
Julien Duval, Paul Gauduchon, Emmanuel Giroux, François Laudenbach,
Manuel Samuelides, Claude Viterbo et, tout spécialement, Yakov
Eliashberg et Jean-Claude Sikorav pour leur travail attentif et
leurs suggestions.
\end{small}

\section{\'{E}noncé et plan}
\label{enonceplan}

Nous cherchons à représenter des classes d'homologie de dimension 2
sur une surface projective $\vX\subset\CPM$ par des surfaces
lagrangiennes lisses. La variété \vX\ est supposée lisse, minimale,
avec $c_{1}(\vX)^{2}>0$~: le fibré canonique $\fcX$ est donc
ample~\cite{BHPV}, et l'on fait en outre l'hypothèse qu'il existe un
diviseur canonique $\Sigmaz$ {\it lisse} et {\it connexe}. \vX\ est
munie de la métrique \gz\ et de la structure symplectique $\omegaz$
induites par la métrique de Fubini-Study sur \CPM~: le triplet
$(\Jz,\gz,\omegaz)$ est compatible (au sens où
$\omegaz(\vectx,\vecty) = \gz(\Jz\vectx,\vecty)~; \hspace{2mm}
\vectx,\vecty\in\Gamma(T\vX)$, $\Jz$ presque-complexe, $\gz$
riemannienne). Enfin, on suppose que la classe $[\omegaz]$ est duale
d'un multiple de $[\Sigmaz]$, ce qui arrive en particulier si
$\omegaz$ est induite par un plongement multicanonique.\\
Pour ne pas surcharger l'écriture, nous noterons simplement
$\|\cdot\|$ les normes pour la métrique \gz.\\
Rappelons que si ${\bf V}$ est un espace vectoriel orienté réel de dimension 4,
$\grass$ la grassmannienne de ses 2-plans orientés et $\mathsf{c}$ une
structure conforme sur ${\bf V}$, alors \grass\ est canoniquement isomorphe à
$\Splus\times\Smoins$, où $\Splus$, respectivement $\Smoins$, sont les sphères
``à l'infini'' dans l'espace des 2-formes auto-duales, respectivement
anti-autoduales pour $\mathsf{c}$. La sphère \Splus\ paramètre les structures
complexes compatibles avec l'orientation et la structure conforme $\mathsf{c}$
sur \espV~; une telle structure complexe $J$ étant fixée, \Smoins\ paramètre
alors les droites $J$-complexes de \espV. Pour tout triplet $(\Jz,\gz,\omegaz)$
cette décomposition en sphères correspond à la décomposition des puissances
extérieures
\[
\Lambda^{+}=\Lambda^{2,0}\oplus\Lambda^{0,2}\oplus\Lambda^{0}.\omegaz~;
\hspace{5mm} \Lambda^{-}=\Lambda^{1,1}_{0}, \] où les éléments de
$\Lambda^{1,1}_{0}$ sont les formes de bidegré $(1,1)$ orthogonales à
$\omegaz$.\\
\indent Considérons une forme différentielle holomorphe $\varphi$ sur \vX, de
type (2,0), qui s'annule sur $\Sigmaz$, et soit $\omegap$ la partie réelle de
$\varphi$.\\
\indent L'idée est la suivante. Dans le complémentaire de $\Sigmaz$, $\omegap$
est une forme symplectique. Elle est auto-duale pour \gz~; sur
$\vX\setminus\Sigmaz$, rempla\c{c}ant \gz\ par la métrique
$\gp=(\|\omegap\|/\sqrt{2}).\gz$, conformément équiva\-lente à \gz, $\omegap$
est de norme $\sqrt{2}$ pour \gp~\cite{Wei}. Ainsi, \gp\ et \omegap\ sont
compatibles, et définissent une structure presque-complexe $J$ sur
$\vX\setminus\Sigmaz$.  De plus, $J$ est une isométrie pour \gz\ et les
structures $J$ et \Jz\ sont orthogonales.\\
\indent Une courbe dans $\vX\setminus\Sigmaz$ pseudo-holomorphe pour $J$ sera
alors lagrangienne pour $\omegaz$. En effet, si $\Sigma$ est une telle courbe,
$J$ envoie $T\Sigma$ sur lui-même~: pour tout $\ptx\in\Sigma$,
$\omegaz(\vectx,J\vectx)$ est nul, d'où le résultat.\\
\indent Si $\gamma$ est une 2-forme anti-auto-duale pour \gz, donc de bidegré
$(1,1)$, orthogonale à $\omegaz$, toute 2-forme
$\omega=\omegap+\varepsilon.\gamma$ avec $\varepsilon$ assez petit sera encore
symplectique hors d'un voisinage de $\Sigmaz$. En général $\omega$ ne sera plus
auto-duale pour \gz~, mais l'opérateur $J$ définira encore une structure
presque-complexe adaptée à $\omega$~:
\begin{lemme}\label{stablejt}
Soit $(\Jz,\gz,\omegaz)$ un triplet compatible sur \vX~; $\omegap$ la partie
réelle d'une $(2,0)$-forme holomorphe (auto-duale) sur $\vX$, avec un zéro sur
le diviseur canonique $\Sigmaz$, $\mbox{$\mathsf
g$\textsf{'}}=\frac{\|\omegap\|}{\sqrt{2}}\gz$ hors de \Sigmaz, et $J$ la
structure presque-complexe associée à $\omegap$ et à $\mbox{$\mathsf
g$\textsf{'}}$ comme précé\-demment.\\
\indent Alors, pour toute 2-forme $\gamma$ anti-auto-du-\\ale pour \gz\ et tout
$\varepsilon>0$ tels que $\|\varepsilon.\gamma\|<\|\omegap\|$ en tout point
hors d'un voisinage $U$ de $\Sigmaz$, $\omega=\omegap+\varepsilon.\gamma$ et
$J$ sont encore compatibles. Elles définissent une métrique $\g$ telle que le
triplet $(J,\g,\omega)$ est compatible. En particulier, $\omega$ est auto-duale
pour la structure conforme associée à $\g$.
\end{lemme}
En effet, la première assertion résulte du fait que $\gamma$ étant
anti-autoduale est de type $(1,1)$ pour $J$ aussi bien que pour $J_0$. La
seconde assertion dit qu'en tout point $\ptx$ hors de $U$, pour tout $\vecx \in
T_{\enx}\vX$, $\omega(\vecx,J\vecx)$ ne s'annule que si $\vecx=0$~; mais
l'hypothèse $\|\varepsilon.\gamma\|<\|\omegap\|$ interdit que
$\omegap(\vecx,J\vecx)=\frac{\|\omegap\|}{\sqrt{2}}\gz(\vecx,\vecx)$ soit
compensé par $\varepsilon.\gamma(\vecx,J\vecx)$, d'où le résultat.\\

\paragraph{} Rappelons que le module d'{\it homologie à distance finie} de
Picard est défini comme l'image $E_{1}$ de
$\HdZ{\vX\setminus\Sigmaz}$ dans $\HdZ{\vX}/\mathrm{Torsion}$.

{\it Remarque~:} A chaque courbe simple $\gamma$ plongée dans
$\Sigmaz$, Picard associe un tore $T_{\gamma}$ bordant un tore plein
constitué de disques holomorphes centrés sur $\gamma$. Ces tores
$T_{\gamma}$ sont lagrangiens. De plus, on sait qu'ils engendrent un
supplémentaire de $E_{1}$ dans le quotient
$\HdZ{\vX}/\mathrm{Torsion}$~\cite{Lef24,Lam81}.

\begin{definitio}
$\tilde{C}^{+}$ désigne l'ensemble des formes
$\lambda\omega=\lambda(\omegap+\varepsilon.\gamma)$,
$\lambda\in\R_{+}^{\ast}$, $\omegap+\varepsilon.\gamma$ comme
ci-dessus, telles qu'il existe un voisinage tubulaire $U$ de
$\Sigmaz$ en dehors duquel $\omega$ est symplectique et adaptée à
$J$. Notons $C^{+}_{\Z}$ (respectivement $C^{+}$) l'ensemble des
classes d'homologie de dimension 2 dans $E_{1}$ (respectivement,
dans $E_{1}\otimes\R$) dont la classe duale a son représentant de
Hodge dans $\tilde{C}^{+}$.
\end{definitio}
Notons que $C^{+}_{\Z}$ est non vide et constitué de classes
d'auto-intersection strictement positive.

\paragraph{} Le principal
résultat obtenu ici concerne le cas où le genre géomé\-trique de
$\vX$ est impair:
\begin{theoreme}\label{theotrois} Soit \vX\ une surface
projective complexe minimale, de genre géométrique impair, avec
$c_{1}(\vX)^{2} > 0$, et telle qu'il existe un diviseur canonique
$\Sigmaz$ sur \vX, lisse, connexe. Supposons \vX\ muni de la
structure symplectique induite par un plongement multicanonique.
Alors tout élément de $C^{+}_{\Z}$ admet un multiple qui se
représente par une surface lagrangienne $\Sigma$ plongée dans
$\vX\setminus\Sigmaz$.
\end{theoreme}

{\it Remarque~:} en vertu du théorème de Bertini (cf.\cite{GrHa}), l'hypothèse
sur $\Sigmaz$ équivaut à l'absence de point de base pour le système linéaire
canonique, ou encore à la régularité de l'application canonique. Nos autres
hypothèses disent que le faisceau canonique est ample, mais il peut ne pas être
très ample.

\paragraph{} Dans leur article {\it Embedded surfaces and
the structure of Donaldson's polynomial invariants}~\cite{KrMr95}, Kronheimer
et Mrowka démontrent que si une variété \vX\ de dimension 4, telle que
$b_{1}(\vX)$ est nulle, $b^{+}(\vX)$ est impair et strictement supérieur à 1,
contient une surface plongée tendue de genre supérieur ou égal à 2, alors \vX\
est de type simple. Ils conjecturent que, dans le cas particulier où \vX\ est
de genre géométrique impair, \vX\ contient une telle surface plongée tendue de
genre supérieur ou égal à 2. Notre théorème résout cette conjecture dans le cas
où \vX\ possède un diviseur canonique lisse et connexe; on peut donc déduire du
théorème de Kronheimer et Mrowka le résultat suivant~:
\begin{corollaire}
Soit \vX\ une surface projective complexe minimale, de genre
géométrique impair, avec $b_{1}(\vX)=0$, $c_{1}(\vX)^{2}
> 0$, et telle qu'il existe un diviseur canonique $\Sigmaz$ sur
\vX, lisse, connexe. Alors \vX\ est une variété simple au sens de
Kronheimer et Mrowka.
\end{corollaire}

{\it Remarque~:} si le genre géométrique de $\vX$ est pair, on
démontre l'existence d'une surface duale d'un multiple de
$[\omega]$, réunion de deux surfaces à bord $\Sigma'$ et $\Sigma''$,
telles que $\Sigma'\subset\vX\setminus U$ est lagrangienne et
$\Sigma''$ est contenue dans $U$. Mais cela résulte aussi des
constructions de Lefschetz. Il serait intéressant de voir, à partir
de notre étude, dans quels cas la surface obtenue pourrait
effectivement être disjointe de $U$ en restant lagrangienne.

\paragraph{}
Pour établir le théorème, nous allons retirer un voisinage tubulaire $U$ de
$\Sigmaz$. Obtenant ainsi une variété \vXz\ à bord $\vY=\partial U$, nous la
complèterons par une variété symplectique $\vW$ de manière à prolonger $\omega$
et $J$ en une forme symplectique et une structure complexe, encore notées
$\omega$ et $J$, sur la variété compacte sans bord $\vXt=\vXz\cup_{\vY}\vW$.

Nous déduirons du théorème de Donaldson~\cite{Don96} l'existence d'une courbe
pseudo-holo\-morphe $\Sigma$ sur $(\vXt,J)$, duale de $k\omega$ pour $k$ assez
grand, mais rencontrant {\it a priori} $\vW$. Restera le point délicat de
disjoindre $\Sigma$ de $\vW$. Nous ferons alors appel aux résultats plus précis
de Donaldson~\cite{Don96,Don99} et d'Auroux~\cite{Aur97,Aur01}, ainsi qu'à une
construction de J.~Duval~\cite{Duv}.

\paragraph{}
Indiquons brièvement le plan que nous allons suivre pour la
démons\-tration.

1) La structure $\omega$ est étudiée au voisinage de $\Sigmaz$~; en
particulier, nous montrons que le bord $\vY$ d'un voisinage
tubulaire $U$ de $\Sigmaz$ est convexe (\cad\ que \vXz\ est {\it
concave}), et nous identifions la structure de contact $F$ induite
par $J$ sur $\vY$ (section~\ref{voisdivcan}).

2) La structure $F$ est une ``structure de contact Spin'', et $\vY$
est une ``fibration Spin'' au-dessus de $\Sigmaz$.
En~\ref{spinatiyah} sont étudiées les structures Spin sur des
surfaces de Riemann, d'après Atiyah. La section~\ref{remplispin}
établit que si la structure Spin sur $\Sigmaz$ induite par $\varphi$
est {\it bordante} (c'est le cas si le genre géométrique $p_{g}$ de
\vX\ est impair), alors $(\vY,F)$ est le bord d'une variété de Stein
$W$, {\it i.e.}, la structure de contact $F$ est holomorphiquement
remplissable.

Sur le complémentaire du 2-squelette $S$ de $W$, nous prolongeons
$\omegaz$ en une forme symplectique orthogonale à $\omega$ . Nous
prouvons que $S$ est réunion de surfaces $\omega$-lagrangiennes
immergées, à croisements normaux.

3) \`{A} la variété symplectique $(\vXt,\omega)$ ainsi obtenue s'applique le
théorème de Donaldson~: soit $(\vX;J,\g,\omega)$ une variété symplectique avec
$[\omega]$ entière, $(J,\g,\omega)$ compatible, $L$ un fibré en droites
complexes sur $\vX$ tel que $c_{1}(L)=[\omega]$~; alors, pour $k$ suffisamment
grand, il existe des sections de $L^{\otimes k}$ transverses à la section nulle
et asymptotiquement holomorphes.

Notons $\Sigma$ le lieu des zéros d'une telle section. Nous
commençons par démontrer que, pour $k$ assez grand, toutes ces
surfaces symplectiques pour $\omega$ peuvent être prises
$\omegaz$-lagrangiennes (en fait, lagrangiennes pour $\omegaz$ sur
$\vXt\setminus S$ et $J$-holomorphes au voisinage de $S$) .

4) Suivant une démonstration par J.~Duval~\cite{Duv} de la convexité
rationnelle des surfaces isotropes des variétés kählériennes, nous arrivons,
pour $k$ grand, à construire une surface plongée $\Sigma$
$\omegaz$-lagran\-gienne, disjointe de $S$, homologue au dual de Poincaré de
$k[\omega]$. Enfin, $\Sigma$ est déformée à l'aide d'un champ de Liouville de
$\omegaz$ pour la rendre plongée dans $\vXz$.

({\it Remarque~:} pour le quatrième point, nous pourrions aussi
faire appel aux méthodes de J.~P.~Mohsen~\cite{Moh,AGM}.)

\section{Au voisinage du diviseur canonique}\label{voisdivcan}

Dans cette section, pour ne pas alourdir les notations, nous
noterons $\Sigma$ la surface $\Sigma_{0}$.

\subsection{$N_{\Sigma}$ racine de \fcS}\label{racinefibrcan}
Remarquons tout de suite que le fibré normal de $\Sigma$ dans \vX\ est une
racine carrée du fibré canonique \fcS. En effet, si $\Sigma$ est déterminée par
une famille d'équations locales $\{\vU_{\alpha},h_{\alpha}\}$, le fibré en
droites associé, $L_{\Sigma}$ (défini par le cocycle
$g_{\alpha\beta}=h_{\alpha}/h_{\beta}\in\holostar(\vU_{\alpha}\cap\vU_{\beta})$),
est isomorphe à \fcX. Les $dh_{\alpha}$ sont des sections locales holomorphes
sans zéros du fibré conormal de $\Sigma$, $N_{\Sigma}^{\ast}$, au-dessus des
ouverts $U_{\alpha}\cap\Sigma \subset \Sigma$ et vérifient
$dh_{\alpha}=g_{\alpha\beta}.dh_{\beta}$~: elles définissent donc une section
globale sans zéro de $N_{\Sigma}^{\ast}\otimes L_{\Sigma}$. Ainsi
$(L_{\Sigma})_{|\Sigma}=N_{\Sigma}$. De cette égalité et de la formule
d'adjonction se déduit le résultat. D'après Atiyah~\cite{Ati71}, ceci signifie
que le fibré normal de $\Sigma$ dans \vX\ définit une structure Spin sur
$\Sigma$.

Rappelons que le genre géométrique $p_{g}=h^{2,0}$ de \vX\ est la
dimension de l'espace des sections holomorphes du fibré canonique
\fcX. Or deux sections possédant le même diviseur sont
proportionnelles, donc la dimension de l'espace des sections
holomorphes de $N_{\Sigma}$ est $p_{g}-1$~; elle est paire lorsque
le genre géométrique de $\vX$ est impair.



\subsection{Structures de contact au bord du voisinage $U$}\label{bordvoisinage}

Soit $\Sigma'$ le germe de $\Sigma$ en un de ses points.

 Dans des coordonnées locales holomorphes $(z,w)$ où $\Sigma'$
s'écrit $\{w=0\}$, soit $\varphi$ une 2-forme méromorphe dont
$\Sigma'$ est un pôle d'ordre $k\geq 2$ ou un zéro d'ordre $k\geq
1$; autrement dit
\[\varphi=w^{k}\,dz\wedge dw~; k\in\Z, k\neq 0, -1.\]
Hors de $\Sigma$, posons $w=\rho e^{it}$ et considérons
$Y_{\epsilon'}$ d'équation $\{\rho=\epsilon'\}$, bord d'un voisinage
tubulaire $U_{\epsilon'}$ de $\Sigma'$. Soit $g'_{0}+i\omegap_{0}$
le produit hermitien sesquilinéaire standard dans ces coordonnées;
notons $\omegap$ la partie réelle de $\varphi$. La forme $\omegap$
et la métrique $g''=\frac{\|\omegap\|}{\sqrt{2}}g'_{0}$ sont
compatibles et définissent donc une structure quasi-complexe $J'$.
Nous voulons établir le lemme suivant~:

\begin{lemme}\label{cylindrelocal}
Les droites complexes pour $J'$ forment une structure de con\-tact
sur $Y_{\epsilon'}$, revêtement à $(k+1)$ feuillets de la structure
des éléments de contact sur $\Sigma'$.
\end{lemme}

Notons $z=x+iy$; une structure presque-complexe $J'$ adaptée à
$\omegap$ hors de $\Sigma'$ est donnée par
\begin{eqnarray*} J'\,dx &=& \rho^{k}\,(\cos(k+1)t \, .\,d\rho -
\sin(k+1)t\, .\,\rho\,dt),
\\ J'\,dy &=& \rho^{k}\,(-\sin(k+1)t\, .\,d\rho -
\cos(k+1)t\, .\,\rho\,dt) \end{eqnarray*} et, par suite,
\begin{eqnarray*}
J'\,d\rho &=& \rho^{-k}\,(-\cos(k+1)t\, .\,dx + \sin(k+1)t\, .\,dy),\\
J'\,dt &=&\rho^{-k-1}\,(\sin(k+1)t\, .\,dx + \cos(k+1)t\, .\,dy).
\end{eqnarray*}

Sur la variété $\vY_{\epsilon'}$, le champ de plans
$F=T\vY_{\epsilon'}\cap J' T\vY_{\epsilon'}$ est défini par
léquation de Pfaff $\alpha=0$ où $\alpha =
-\frac{\epsilon'^{k+1}}{k+1}\, [\cos(k+1)t\,dx - \sin(k+1)t\,dy]$.
Comme $\alpha\wedge d\alpha=
\frac{\epsilon'^{2(k+1)}}{(k+1)^{2}}\,dx\wedge dy\wedge dt$, $F$ est
une structure de contact sur $\vY_{\epsilon'}$ qui oriente
$\vY_{\epsilon'}$ comme bord de $U_{\epsilon'}$~; comme, de plus,
$d\alpha|_{T\vY_{\epsilon'}}\equiv\omegap|_{T\vY_{\epsilon'}}$,
$\alpha$ est adaptée à $\vY_{\epsilon'}$ et à $\omegap$, et
$\vY_{\epsilon'}$, bord de $U_{\epsilon'}$, est du type de contact
dans \vX.

Pour tout $k\neq -1$, la normale sortante \vectn\ de
$\vY_{\epsilon'}$ étant $\partial_{\rho}$ et le champ de Reeb
\vectxi\ de $\alpha$ étant
$-\frac{k+1}{\epsilon'^{k+1}}\,[\cos(k+1)t\,\partial_{x} -
\sin(k+1)t\,\partial_{y}]$, on calcule
\[\omegap(\vectn,\vectxi) = \frac{k+1}{\epsilon'}~.\] $(\vY_{\epsilon'},F)$
est concave par rapport à $\omegap$ pour tout entier $k\leq -2$, et
convexe par rapport à $\omegap$ pour tout entier $k\geq 0$.
\qed\newline

\paragraph{}
L'application que nous ferons de ce lemme ne concerne que le cas $k=1$. Même
localement, on ne peut supposer que $g_{0}$ co\"{\i}ncide avec une métrique
plate comme l'est $g'_{0}$; par contre une telle égalité est vraie à l'ordre 2
en un point donné. On
 en déduit le lemme suivant~:

\begin{lemme}\label{contactrevetement}
Dans un voisinage tubulaire $U$ de $\Sigma$, soit $\varrho$ la
distance à $\Sigma$ et $Y_{\epsilon}$ la sous-variété d'équation
$\varrho=\epsilon$. Les droites $J$-complexes forment sur
$Y_{\epsilon}$ un champ de plans de contact $F_{\epsilon}$ isomorphe
à un revêtement double du champ de plans de la structure standard de
$ST\Sigma$.
\end{lemme}

Fixons une fois pour toutes un revêtement du fibré en cercles unités
pour une métrique sesquilinéaire  $g'_{0}$ sur le fibré normal
$N_{\Sigma}$~:
\[ (Y \subset N_{\Sigma}) \stackrel{\pi}{\longrightarrow} (ST\Sigma
\subset T\Sigma)\] et notons $F'$ la structure sur $Y$ qui s'en
déduit.

Dans chaque ouvert trivialisant $U^{(\alpha)}$ du fibré $ST\Sigma$,
le lemme~\ref{cylindrelocal} fournit une structure $J'^{(\alpha)}$
telle que $F'$ s'identifie à la structure
$F'\,_{\epsilon'}^{(\alpha)}$ du lemme sur le bord
$Y_{\epsilon'}^{(\alpha)}$ d'un voisinage
$U_{\epsilon'}^{(\alpha)}$.

Par ailleurs, le flot géodésique pris aux temps $\epsilon$
successifs fournit une famille à un paramètre, analytique réelle, de
difféomorphismes $\varphi_{\epsilon}$ de $Y$ sur $Y_{\epsilon}$. La
structure $\varphi_{\epsilon}^{\ast}F_{\epsilon}$ restreinte à
$U^{(\alpha)}$ est $C^{1}$-proche de $F'\,_{\epsilon'}^{(\alpha)}$.
Grâce au théorème de Moser et Gray ({\it cf.}~\cite{Ben83}), les
structures $\varphi_{\epsilon}^{\ast}F_{\epsilon}$ et $F'$ sont
alors difféomorphes. \qed\newline

Toujours en vertu du théorème de Moser et Gray, les lemmes
\ref{cylindrelocal} et \ref{contactrevetement} restent valables si
$\omegap$ est remplacée par une forme symplectique
$\omega=\omegap+\varepsilon.\gamma$, où $\gamma$ est une (1,1)-forme
anti-auto-duale (orthogonale à $\omegaz$) et $\varepsilon$ est
suffisamment petit.


\section{Structures Spin sur les surfaces de Riemann}\label{spinatiyah}


Rappelons brièvement la définition d'une structure Spin. Soit \vV\
une variété réelle riemannienne orientée de dimension $n\geq 2$, et
$P_{SO}(\vV )$ le fibré [principal à droite, de groupe structural
$SO(n)$] de ses repères orthonormés. (Pour $n=2$, nous convenons que
$Spin(2)$ est $SO(2)$ vu comme revêtement double connexe de
$SO(2)$.)

Une {\it structure Spin} sur \vV\ est la donnée d'un $Spin(n)$-fibré
principal $\dot{P}$ et d'un revête\-ment double $\dot{\xi}$,
morphisme de fibrés $Spin(n)$-équivariant~:
  \[ \dot{P} \stackrel{\dot{\xi}}{\longrightarrow} P_{SO}(\vV ). \]

Une condition nécessaire et suffisante pour que \vV\ admette une
structure Spin est la nullité de sa deuxième classe de
Stiefel-Whitney $w_{2}(\vV)$. On dit que deux structures Spin sont
équivalentes s'il existe un diagramme commutatif d'isomorphismes des
fibrés principaux correspondants.

Si $w_{2}(\vV)=0$, la suite exacte
\[ 0\rightarrow
H^{1}(\vV;\Zsd)\stackrel{\pi^{\ast}}{\rightarrow}H^{1}(P_{SO}(\vV);\Zsd)\rightarrow
H^{1}(SO(n);\Zsd)\stackrel{w}{\rightarrow}H^{2}(\vV;\Zsd), \]
extraite de la suite spectrale associée à la fibration
\[ SO(n)\rightarrow P_{SO}(\vV)\stackrel{\pi}{\rightarrow}\vV, \]
permet de conclure que les structures Spin sur \vV\ sont en
bijection avec les éléments de $H^{1}(\vV;\Zsd)$ (de fait, $w$
envoie le générateur de $H^{1}(SO(n);\Zsd)\cong\Zsd$ sur
$w_{2}(\vV)$ dans la suite spectrale ci-dessus).

Si \vV\ est de dimension paire $2m$, munie d'une structure
presque-complexe, la réduction de $SO(2m)$ à $U(m)$ permet de
montrer, à partir des suites exactes
\[ 0\rightarrow
H^{1}(\vV;\Zsd)\stackrel{\pi^{\ast}}{\rightarrow}H^{1}(P_{U}(\vV);\Zsd)\rightarrow
H^{1}(U(m);\Zsd)\rightarrow H^{2}(\vV;\Zsd) \] et \[ 0\rightarrow
H^{1}(\vV;\Zsd)\stackrel{\pi^{\ast}}{\rightarrow}H^{1}(\det
P_{U}(\vV);\Zsd)\rightarrow H^{1}(U(1);\Zsd)\rightarrow
H^{2}(\vV;\Zsd), \] que les structures Spin sur \vV\ sont en
bijection avec les revêtements doubles du $U(1)$-fibré déterminant
des repères dont la restriction à chaque fibre est $U(1)\rightarrow
U(1), \eit\mapsto e^{2i\theta}$. Par conséquent, si \vV\ est
analytique complexe, l'homomorphisme $\pi^{\ast}$ réalise une
bijection entre les classes de structures Spin sur \vV\ et les
classes d'isomorphisme des couples $(L,\kappa)$, où $L$ est un fibré
en droites holomorphe et $\kappa$ un isomorphisme holomorphe de
$L\otimes L$ dans le fibré canonique $\mathcal{K}_{V}$ de \vV. Si,
de plus, \vV\ est compacte, la structure holomorphe de $L$ détermine
$\kappa$ à une constante multiplicative près, de sorte que les
classes de structures Spin sur une variété complexe compacte \vV\
sont en bijection naturelle avec les classes d'isomorphisme de
fibrés en droites holomorphes $L$ tels que $L\otimes L$ est
isomorphe au fibré canonique $\mathcal{K}_{V}$ de \vV.

Lorsque $m=1$ et que \vV\ est une surface de Riemann compacte
$\Sigma$ de genre $g$, il existe exactement $2^{2g}$ structures Spin
non équivalentes. Le fibré canonique \fcS\ de $\Sigma$ est son fibré
cotangent holomorphe. La caractéristique d'Euler de $\Sigma$,
$2-2g$, étant paire, \fcS\ admet une racine carrée, \cad\ un fibré
en droites holomorphe $L$ tel que $L\otimes L\cong\fcS$. Notons
$\mathcal{R}(\Sigma)$ l'ensemble des classes d'isomorphisme de tels
fibrés, racines carrées de \fcS. L'on définit la fonction
\[ \varphi~: \mathcal{R}(\Sigma)\rightarrow\Zsd,\hspace{5mm}
L\mapsto\dim H^{0}(L)\hspace{3mm}[mod\ 2]. \]

Atiyah démontre la proposition suivante~: \textit{La fonction $\varphi$ est
qua\-dratique et sa forme bilinéaire associée $B_{\varphi}$ s'identifie au
cup-produit sur $H^{1}(\Sigma;\Zsd)$.} On appelle \textit{paire},
respectivement \textit{impaire}, une classe de structure Spin $\spbd$ telle que
$\varphi(\spbd)=0$, resp. $\varphi(\spbd)=1$.

\paragraph{}
Si \vV\ est une variété compacte à bord de dimension $n>1$, une
structure Spin sur \vV\ induit une structure Spin sur
$\partial\vV$~\cite{Mil63}. Une variété compacte $\vV$ munie d'une
structure Spin est {\it Spin-cobordante à 0} (ou, en abrégé, {\it
bordante}) si elle est différentiablement équivalente au bord
$\partial\vV'$ d'une variété Spin, muni de la structure Spin
induite.

Le théorème de l'indice d'Atiyah-Singer ({\it cf.}~\cite{AtSi5},
théorème (3.3), et \cite{Ati71}) permet alors de démontrer le
théorème suivant~: soit $\spbd$ une structure Spin sur $\Sigma$~;
alors $(\Sigma, \spbd)$ est Spin-cobordante à 0 si, et seulement si
$\spbd$ est paire.

Atiyah en déduit le théorème~:

\textit{Une surface de Riemann $\Sigma$ de genre $g$ admet
exactement $2^{g-1}(2^{g}+1)$ structures Spin cobordantes à zéro.}
\newline

{\it Exemples~:}\label{spintore} Pour $\Sigma=T^{2}$, le fibré en
cercles unitaires du cotangent, $ST^{\ast}T^{2}$, s'identifie à
$T^{3}$. Les structures Spin correspondent aux sous-groupes d'indice
2 dans $\pi_{1}(T^{3})=\Z^{3}$ qui s'envoient surjectivement sur
$\pi_{1}(T^{2})$~: ce sont les quatre sous-groupes engendrés par les
triplets
\[ (2e_{0},e_{1},e_{2})~; (e_{0},e_{1}+2e_{0},e_{2})~; (e_{0},e_{1},e_{2}+e_{0})~;
(e_{0},e_{1}+2e_{0},e_{2}+2e_{0}).\] Les trois premières structures
Spin sont cobordantes à zéro, pas la quatrième.

Toujours sur $T^{3}$, avec des coordonnées angulaires
$(\theta,x_1,x_2)$, pour tout $n\geq 1$, notons $\zeta_{n}$ la
structure de contact définie par la 1-forme $(\cos
n\theta)dx_{1}+(\sin n\theta)dx_{2}$. Cette structure est tendue
pour tout $n$. Les trois structures Spin bordantes donnent
$\zeta_{1}$~; la structure non-bordante, $\zeta_{2}$. Eliashberg a
démontré que $\zeta_{2}$ n'était pas holomorphiquement remplissable
(elle ne borde pas une variété de Stein), alors que $\zeta_{1}$ est
la structure standard sur $ST^{\ast}T^{2}$, bord de Stein.
Néanmoins, d'après Giroux ({\it cf.}~\cite{Gir}), $\zeta_{2}$ (de
même que les autres $\zeta_{n}$) borde une variété symplectique
(compacte à bord) du côté convexe.

\section{Remplissage des structures de contact
Spin}\label{remplispin}

L'objet de cette section est de démontrer le résultat suivant~:
\begin{propositio}\label{remplie}
Soit $\spbd$ une structure Spin sur une surface de Riemann $\Sigma$,
cobordante à zéro. Alors la structure de contact correspondante est
holomorphiquement remplissable, c'est-à-dire qu'elle borde une
variété de Stein.
\end{propositio}

\subsection{Spin-cobordisme et genre géométrique}\label{spincobord}

En ce qui concerne notre problème de départ, la condition néces-\\saire et
suffisante pour que la structure Spin induite par $\omega$ sur le diviseur
$\Sigma$ soit cobordante à zéro est que la dimension de l'espace des sections
holomorphes du fibré normal $N_{\Sigma}$ soit paire. Or cette dimension est
égale au genre géométrique diminué de 1, donc la structure Spin considérée
borde si, et seulement si $p_{g}$ est impair. Sous cette hypothèse, nous
obtiendrons donc la variété symplectique $\tilde{X}$ annoncée en introduction
comme conséquence de la proposition \ref{remplie}.

\subsection{Le modèle global du remplissage}\label{modeleglobal}

Atiyah montre que si une structure Spin sur la surface de Riemann
$\Sigma$ est cobordante à zéro, elle borde dans une anse pleine
(\cite{Ati71}, p. 58).
\paragraph{}
 Cette hypothèse étant supposée vérifiée,
donnons-nous à présent une anse pleine $V$ de bord $\Sigma$~; une
structure Spin {\boldmath $s'$} sur $V$ de bord \spbd~; et une
fonction de Morse $f: V\rightarrow\R$ possédant $g$ points critiques
$p_{1},\ldots,p_{g}$ d'indice 1 (au voisinage d'un tel point, $\grad
f$ est rentrant dans 2 directions et sortant dans la troisième), et
un unique point critique $p_{0}$ d'indice 0 (au voisinage duquel le
gradient est rentrant dans toutes les directions). Nous noterons $R$
le fibré des repères orthonormés directs sur \vV\ et $P$ le
revêtement spinoriel associé à $\spbd'$~: $P\rightarrow R$.

Les paragraphes~\ref{modele1}, \ref{modele2} et \ref{modele3}
donneront un fibré en cercles $\tilde{M}$ au-dessus de
$V^{\times}=V\setminus (B_{1}\cup\ldots\cup B_{g}\cup B_{0})$,
plongé dans $P$, et une manière d'attacher à $\tilde{M}$ des boules
$\tilde{B}_{1},\ldots,\tilde{B}_{g},\tilde{B}_{0}$, modèles locaux
de remplissage symplectique, permettant d'obtenir une variété
symplectique $(W,\omega)$, de bord \vY, et une fonction
$F~:W\rightarrow\R$ relevant $f$ au-dessus de $V^{\times}$, fonction
de Liouville pour $\omega$, qui possède un unique point critique
$\tilde{p}_{i}$ dans chaque $\tilde{B}_{i}$.

\subsection{Extension régulière et d'indice zéro} \label{modele1}

 Entre deux
niveaux critiques de $f$, nous remplissons la structure de contact
Spin au bord de la manière suivante. Soient $a$ et $b$ deux valeurs
critiques consécutives de $f$ ($0<a<b$, $]a,b[$ ne contenant pas de
valeur critique)~: $f^{-1}(]a,b[)\subset\vV$ s'identifie au produit
d'une surface riemannienne réelle $S$ par $]a,b[$~; posons
$M_{a,b}:=ST^{\ast}S\times ]a,b[$.

En chaque point de $S\times ]a,b[$, la donnée d'un vecteur tangent à
$S$, de longueur 1, jointe à la direction de $t$ dans $]a,b[$,
définit un point du fibré $R$. Nous avons donc un plongement
canonique $j$ de $M_{a,b}$ dans $R$, et nous notons
$\tilde{M}_{a,b}$ le relèvement de $j(M_{a,b})$ dans l'espace total
du revêtement Spin $P\rightarrow R$.

Le plongement de $M_{a,b}$ dans $T^{\ast}S$ qui envoie $(x,v;t)$ sur
$(x;e^{t}.v)$ munit $M_{a,b}$ d'une structure symplectique
$\varpi$~; soit $\tilde{\varpi}$ la structure qui s'en déduit sur
$\tilde{M}_{a,b}$. D'où l'énoncé suivant~:

\begin{lemme} \label{spinrevet} Pour tout $t\in
]a,b[$, la structure de contact Spin au-dessus de $S\times\{t\}$
borde la structure $\tilde{\varpi}$ sur $\tilde{M}_{a,b}$ au-dessus
d'un produit $S\times [t-\epsilon,t[$, $\epsilon>0$.
\end{lemme}

Nous verrons en \ref{orthoomegaomegaz} comment raccorder ces
structures entre elles loin des points critiques, grâce à des formes
de Liouville.\\

Pour compléter la construction, il faut disposer de modèles locaux
permettant de remplir les structures de contact tracées au-dessus du
bord des $B_{i}$, voisinages des points critiques $c_{i}$.

(Comme $f$ cro\^{\i}t dans la direction du remplissage, la dimension
des variétés dilatantes en un point-selle est 1~: on dira qu'on a
affaire à un point d'indice 1~; un puits pour le gradient est, quant
à lui, d'indice 0.)

Pour le point critique $p_{0}$ de $f$ d'indice 0, choisissons une
petite boule centrée en $p_{0}$ dans $V$, dont le bord est un niveau
de $f$~; au-dessus de ce bord, le lemme~\ref{spinrevet} fournit une
sphère $S^{3}$ dans $P$, munie de la structure de contact standard,
revêtement double de $ST^{\ast}S^{2}$ munie de la structure des
éléments de contact. Cette structure de contact est remplissable de
façon standard.

Reste donc à remplir par des structures symplectiques sur des boules
de dimension 4 au-dessus du voisinage des points critiques d'indice
$1$. Le remplissage de Stein résultera alors des théorèmes
d'Eliashberg~\cite{Eli90a} et de Gompf~\cite{Gom98}.

\subsection{Préparation d'indice 1} \label{modele2}

Soit $p_{i}, i\geq 1$ un tel point critique d'indice 1. Au-dessus
d'un voisinage de $\partial B_{i}$, le champ de Liouville est
projetable sur un champ de vecteurs rentrant dans la boule $B_{i}$
dans deux directions et sortant dans la troisième, puisqu'il relève
$\grad f$ par construction.

D'autre part, le bord de $\tilde{M}$ au-dessus de $\partial B_{i}$
est une variété $Z$ de dimension 3 fibrant en cercles au-dessus
d'une sphère de dimension 2~; notons $\mathcal{F}$ cette fibration
et appelons $\mathcal{P}_{1}$ le champ de plans tangents dans $V$
défini le long de $\partial B_{i}$ par les perpendiculaires au
gradient de $f$.

Identifions $\partial B_{i}$ à la sphère $S^{2}$ de centre $0$ et de
rayon $1$ dans $\mathbb{R}^{3}$ et considérons le champ
$\mathcal{P}_{0}$ de plans tangents à la sphère $S^{2}$ orientés par
les normales sortantes. La latitude des points de contact
correspondants étant repérée par l'angle $\theta$ (égal à $0$ au
pôle sud et à $\pi$ au pôle nord), nous appliquons à chaque plan une
rotation $\varrho_{t}$ d'angle $2\theta$ autour de la droite passant
par le point de contact et parallèle à l'axe des pôles ($t\in
[0,1]$). Le champ $\mathcal{P}_{1}$ est alors homotope à
$\varrho_{1}(\mathcal{P}_{0})$, pour lequel les plans aux pôles sont
orientés par les normales sortantes, les plans le long de l'équateur
par les normales {\it rentrantes}.

Relevant les deux champs de plans homotopes $\mathcal{P}_{0}$ et
$\mathcal{P}_{1}$ dans $S^{3}$ et dans $Z$ respectivement, nous
voyons que la fibration en cercles $(Z,\mathcal{F})$ est homotope,
donc isotope à la fibration de Hopf de $S^{3}$.

Ces observations vont nous permettre de compléter le remplissage
symplectique à l'aide du modèle local suivant.

\subsection{Modèle local pour le remplissage symplectique au
voisinage des points critiques d'indice $1$} \label{modele3}

Soient $z=x+iy$ et $w=u+iv$ les coordonnées standard sur \Cd, $J$ la
structure presque-complexe donnée par la multiplication par $i$.

\begin{lemme}
Il existe, au voisinage de l'origine dans $\Cd$, une fonction
strictement pluri-sous-harmonique $\varphi$ possédant la propriété
suivante~: si $\hat{\omega}$ est la forme symplectique
$dd^{J}\varphi$, le champ de Liouville pour $\hat{\omega}$ le long
d'une sphère $S^{3}$ centrée en 0, de rayon $r$ suffisamment petit,
est projetable suivant la fibration de Hopf et se projette sur la
restriction à $S^{2}$ dans $\R^{3}$ du champ standard pour un point
critique d'indice 1.
\end{lemme}
La projection $h$ de $S^{3}$ sur $S^{2}$ étant donnée
par \begin{eqnarray*} && (x+iy,u+iv) \\
&& \mapsto
(x_{1},x_{2},x_{3})=(2(xy-uv),2(xv+uy),(x^{2}+u^{2})-(y^{2}+v^{2})),
\end{eqnarray*} le tore $T$ s'envoie sur le cercle $r^{2}S^{2}\cap\{x_{3}=0\}$
de rayon $r^{2}$, les cercles $C_{1}$ et $C_{2}$ s'envoyant
respectivement sur les points $(0,0,r^{2})$ et $(0,0,-r^{2})$.
(L'application $h$ est une fibration de Hopf orthogonale à la
fibration donnée par la structure standard de $\Cd$.)

Nous allons chercher à construire $\varphi$ sous la forme
$\varphi(z,w)=F(x^{2}+u^{2},y^{2}+v^{2})$. L'idée est de prendre un
lissage du carré de la distance à $\{x=u=0\}\cup\{y=v=0\}$, modifié
par une petite perturbation qui permette d'obtenir les propriétés
voulues pour le champ $\grad\varphi$ (pour la métrique standard de
\Cd). Plus précisément, soit $g:\R\rightarrow [0,1]$ une fonction de
classe \Cinf, qui vaut $0$ sur $]-\infty,-\eta]$ ($\eta>0$ petit),
$1$ sur $[\eta,+\infty[$, et est strictement croissante sur
$[-\eta,\eta]$, de dérivées première et seconde bornées, prenant la
valeur $1/2$ en $0$~; et soit $\epsilon>0$ petit. Définissons la
fonction $\varphi$ sur \Cd\ par
\begin{eqnarray*}
&& \varphi(x+iy,u+iv)=
\\ && \hspace{2cm}
(x^{2}+u^{2}).[1-g(x^{2}+u^{2}-y^{2}-v^{2})-\epsilon]
\\ && \hspace{4cm} +(y^{2}+v^{2}).[g(x^{2}+u^{2}-y^{2}-v^{2})-\epsilon].
\end{eqnarray*} Alors~:
\begin{eqnarray*}
\hat{\omega}&=&dd^{J}\varphi
\\&=& 2.(1-2\epsilon)\,(dx\wedge dy+du\wedge dv)\\
&&  \hspace{1cm} +
  4.[(x^{2}+u^{2}-y^{2}-v^{2}).g''(x^{2}+u^{2}-y^{2}-v^{2})
  \\ &&  \hspace{2cm} +2.g'(x^{2}+u^{2}-y^{2}-v^{2})]
  \\
&&  \hspace{3cm} .[-(x^{2}+y^{2})\,dx\wedge
dy-(u^{2}+v^{2})\,du\wedge
  dv \\ && \hspace{35mm}
-(xu+yv)\,(dx\wedge dv-dy\wedge du)
\\ && \hspace{4cm}
+(yu-xv)\,(dx\wedge du+dy\wedge
  dv)]~;
\end{eqnarray*}
on constate que pour $\epsilon<1/2$, $x^{2}+y^{2}+u^{2}+v^{2}<r^{2}$
assez petit et $\eta\ll r^{2}$, la forme $\hat{\omega}$ est proche
de $2\,(dx\wedge dy+du\wedge dv)$, elle reste donc symplectique et
la formule
$(\vecv,\vecw)\mapsto\hat{\omega}(\vecv,J\vecw)+i\hat{\omega}(\vecv,\vecw)$
définit une structure hermitienne définie positive. Par suite
$\varphi$ est strictement pluri-sous-harmonique au voisinage de $0$.

D'autre part, le champ $\grad\varphi$ se projette bien comme
spécifié par le lemme. En effet, notant toujours
$(x_{1},x_{2},x_{3})=(2(xy-uv),2(xv+uy),(x^{2}+u^{2})-(y^{2}+v^{2}))$~:
\begin{eqnarray*}
-\grad\varphi & = &
2x.(x_{3}g'(x_{3})-1+g(x_{3})+\epsilon)\:\partial_{x} \\
& & +2y.(-x_{3}g'(x_{3})-g(x_{3})+\epsilon)\:\partial_{y} \\
& & +2u.(x_{3}g'(x_{3})-1+g(x_{3})+\epsilon)\:\partial_{u} \\
& & +2v.(-x_{3}g'(x_{3})-g(x_{3})+\epsilon)\:\partial_{v}~;
\end{eqnarray*}
si $\pi:\R^{4}\setminus\{0\} \rightarrow \R^{3}\setminus\{0\}$ est
la projection définie par $\pi(\rho.(z,w))=\rho^{2}.h(z,w)$
($\rho\in\R_{+}^{\ast}$, $(z,w)\in S^{3}$), alors $-\grad\varphi$
est projeté par $\pi$ sur le champ de vecteurs
\begin{eqnarray*}
d\pi.-\grad\varphi & = & 2(1-2\epsilon)x_{1}\,\partial_{x_{1}} +
2(1-2\epsilon)x_{2}\,\partial_{x_{2}} \\
& & - 4(\epsilon
x_{3}+\rho^{2}x_{3}g'(x_{3})+\rho^{2}g(x_{3})-(x^{2}+u^{2}))\,\partial_{x_{3}},
\end{eqnarray*}

lequel est conjugué au gradient dans $\Rt$ de la fonction
\[(x_{1},x_{2},x_{3})\mapsto -x_{1}^{2}-x_{2}^{2}+x_{3}^{2}.\] \qed

\subsection{Fin du remplissage. Orthogonalité des formes $\omega$ et $\omegaz$ hors d'un
squelette $\omega$-lagrangien}\label{orthoomegaomegaz}

Revenant au remplissage global, nous pouvons supposer que tous les
$p_{i}$, $1\leq i\leq g$, sont au même niveau $c_{1}>0$ pour $f$, et
que $f(p_{0})=c_{0}>c_{1}$. (0 est le niveau de $\Sigma$.)
Choisissons de petites boules $B_{1},\ldots,B_{g},B_{0}$ autour de
ces points.

Le point $\tilde{p}_{0}$ est un point critique d'indice 0. En
revanche, les $\tilde{p}_{i}$, $i\geq 1$, ne sont pas des points
critiques non-dégénérés. Nous pourrions modifier $F$ en une fonction
de Morse pluri-sous-harmonique ({\it cf.} Biran~\cite{Bir01}) pour
décrire un complexe simplicial isotrope dans $W$~; mais il est plus
simple de conserver $F$ en l'état~: près de $\tilde{p}_{i}$ ($i\geq
1$), il existe deux disques lagrangiens $\Delta_{1}^{(i)}$ et
$\Delta_{2}^{(i)}$ contenant $\tilde{p}_{i}$, transversalement
attractants pour $\grad F$, préservés par $\grad F$, lequel est
dilatant sur ces disques. En suivant les lignes de flot de $\grad
F$, les bords $C_{1}^{(1)}, C_{2}^{(1)}, \ldots, C_{1}^{(g)},
C_{2}^{(g)}$ de ces disques se retrouvent tous au bord de
$\tilde{B}_{0}$ avant de s'évanouir en $\tilde{p}_{0}$. Une
perturbation générique de $F$ assure que les disques d'écrasement en
$\tilde{p}_{0}$ de ces cercles soient deux à deux transverses. (La
figure ci-dessus représente le cas $g=2$.)



\begin{figure}
\begin{center}
\resizebox{0.75\textwidth}{!}{\includegraphics{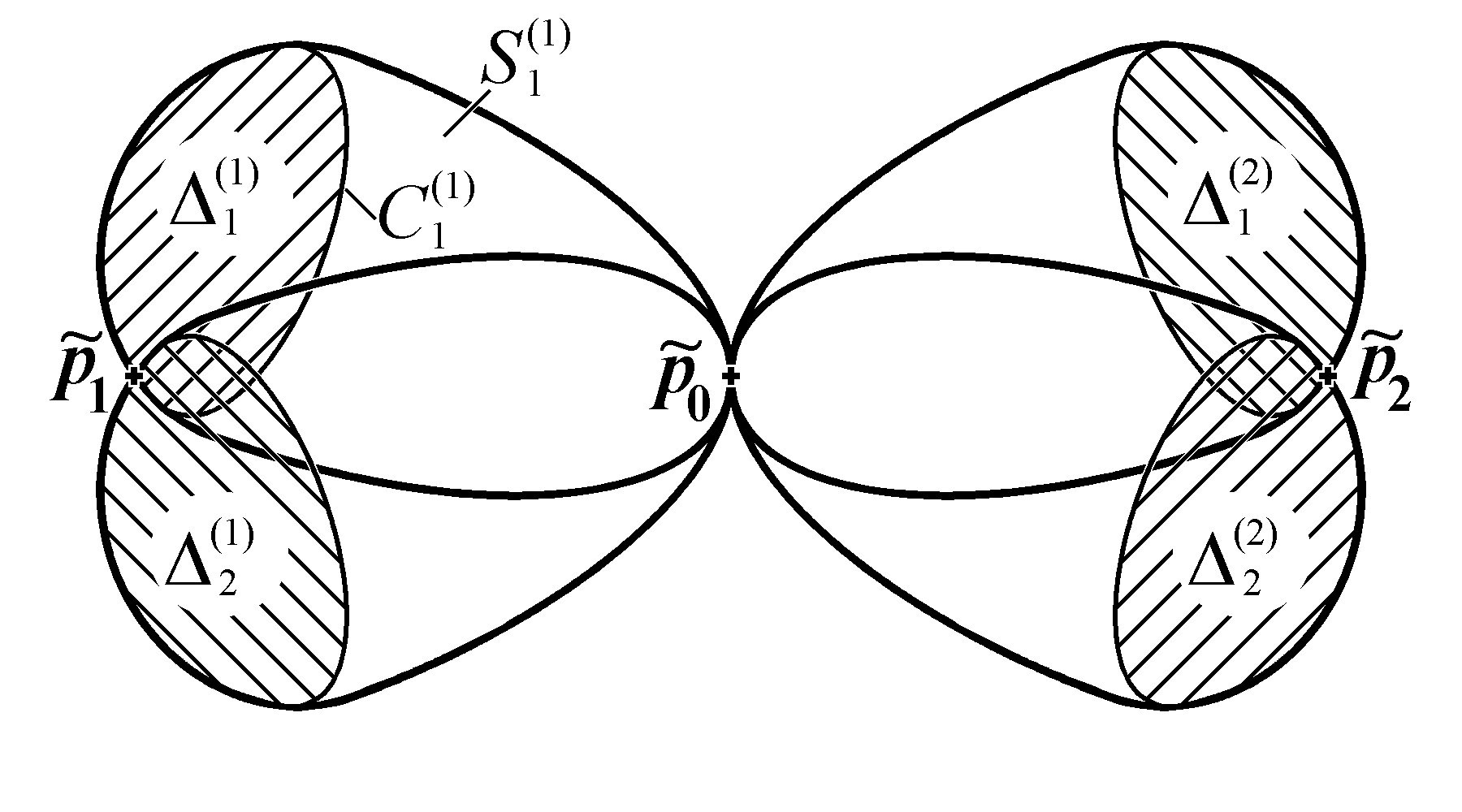}}
\caption{Squelette de $W$} \label{fig:SqueletteW}
\end{center}
\end{figure}

D'où le résultat suivant~:
\begin{lemme}
Le remplissage $W$ est rétracté par un champ de Liouville sur un
bouquet en $\tilde{p}_{0}$ de $g$ couples de sphères
$S_{1}^{(i)}\cup S_{2}^{(i)}$, lagrangiennes (pour $\omega$),
s'intersectant transversalement en $\tilde{p}_{0}$ et en un unique
autre point ($\tilde{p}_{i}$).
\end{lemme}

\begin{definitio}
Le \textup{squelette $S$ de $W$} est la réunion des sphères
$S_{k}^{(i)}$ ($1\leq k\leq 2, 1\leq i\leq g$). C'est une
sous-variété immergée de $W$, $\omega$-lagrangienne, à croisements
normaux.
\end{definitio}

\paragraph{} En dehors de $S$, les lignes de gradient de $F$ permettent de
prolonger $\omegaz$ en une structure symplectique qui reste
orthogonale à $\omega$.



\begin{figure}
\begin{center}
\resizebox{0.5\textwidth}{!}{\includegraphics{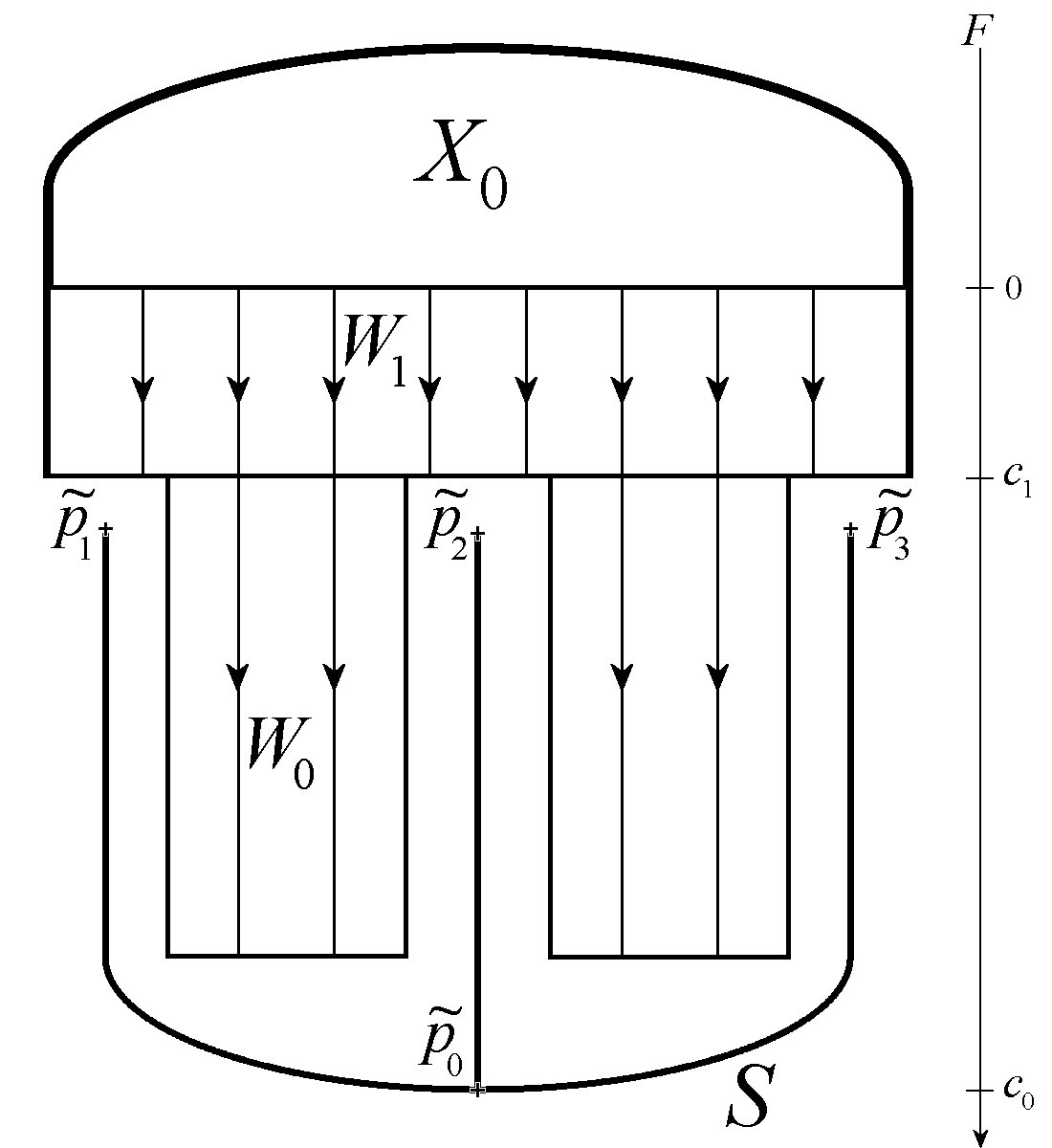}}
\caption{Extension de $\omegaz$ à $\vXt\setminus S$}
\label{fig:extomega}
\end{center}
\end{figure}

Entre les niveaux 0 et $c_{1}-\epsilon$ de $F$ ($\epsilon>0$ petit),
nous avons dans $\vW$ un collier $\vW_{1}=Y\times
[0,c_{1}-\epsilon]$. Le gradient de $F$, $\vectxi=\grad F$, est un
champ de Liouville contractant pour $\omega$~; notons $\lambda$ la
forme de Liouville $\vectxi\prodint\omega$.

Le long de $Y$, il existe une forme exacte $\lambdaz$, de Liouville
pour $\omegaz$, et un champ de Liouville $\vectxi_{0}$ (tel que
$\lambdaz=\vectxi_{0}\prodint\omegaz$) multiple de $\vectxi$~;
$\lambda_{0\;|Y}$ définit alors un plan orthogonal à celui que
définit $\lambda$. Prolongeons canoniquement $\lambdaz$ sur $W_{1}$
de telle sorte que $\lambdaz$ continue à s'annuler à la fois sur
$\vectxi$ et sur un plan orthogonal à celui que définit $\lambda$ et
à la fibration en cercles $W\rightarrow V$~; $\omegaz$ se prolonge
alors sur $W_{1}$ par $\omegaz=d\lambdaz$.

Le champ $\vectxi$, contractant pour $\omega$, est {\it dilatant}
pour $\omegaz$. La même construction peut aussi être poursuivie
entre les niveaux $c_{1}-\epsilon$ et $c_{0}-\epsilon$ de $F$ sur
une partie $W_{0}$ de $W$ qui est saturée par les lignes de gradient
de $F$ évitant les boules $\tilde{B}_{1},\ldots,\tilde{B}_{g}$~;
$W_{0}$ est le complémentaire d'un voisinage de $S$ dans
$F^{-1}([c_{1}-\epsilon,c_{0}-\epsilon])$.

Nous pouvons ainsi définir une forme $\omegaz$ sur $\vXt\setminus
S$~; il est certainement impossible de la prolonger à $S$. Par
contre, la structure $J$, comme $\omega$, se prolonge bien à $\vXt$,
ce qui permettra de parler ultérieure\-ment de variétés
$\omegaz$-lagrangiennes (là où $\omegaz$ est défini) et dont le plan
tangent est $J$-complexe le long de $S$.

{\it Remarque~:} à partir d'un germe de complexification de la
surface singulière $S$, la construction peut être réalisée de sorte
que $J$ soit intégrable au voisinage de $S$, ce que nous supposerons
dorénavant.

\paragraph{} Nous aurons besoin de la précision suivante sur le
complexe $S$. Introduisons d'abord une définition~:
\begin{definitio} \label{poselem}
Nous dirons que $N$ plans lagrangiens vectoriels $\Lambda_{1}$,
\ldots, $\Lambda_{N}$ de $\Cd$ muni de sa structure symplectique
standard $\omegaz$ sont en \upshape{position élémen\-taire} s'il
existe des coordonnées unitaires $(z,w)$ telles que pour
$j=1,\ldots,N$, l'équation de $\Lambda_{j}$ soit
\[ w=e^{\frac{2i(j-1)\pi}{N}}\bar{z}. \]
\end{definitio}

\begin{lemme}
Le complexe $S$ peut être choisi de sorte qu'il existe des
coordonnées $J$-complexes en $\tilde{p}_{0}$ qui mettent les plans
tangents en $\tilde{p}_{0}$ à $S$ ({\it i.e.,} aux sphères
$S^{k}_{(i)}$) en position élémentaire.
\end{lemme}
Nous pouvons choisir la fonction $f$ pour que, en $p_{0}$, les $2g$
trajectoires instables des points $p_{1},\ldots,p_{g}$ aboutissent
dans un même plan tangent de $T_{p_{0}}V$. Dans la sphère $S_{(0)}$
de dimension 3 autour de $\tilde{p}_{0}$ dans $W$, les cercles
d'intersection des $S^{k}_{(i)}$ avec $S_{(0)}$, qui font partie
d'une même fibration de Hopf, se projettent sur un même cercle de
$S^{2}$, et peuvent donc être déformés pour être mis en position
élémentaire. \qed

{\it Remarque~:} en un point $\tilde{p}_{i}$ ($i\geq 1$), il n'y a
que deux plans lagrangiens transverses, tangents à $S$~; il est
toujours possible de choisir des coordonnées symplectiques dans
lesquelles ils aient pour équations $w=\bar{z}$ et $w=-\bar{z}$.

Si $J$ est donnée, $S$ peut aussi être déformée pour que tous ces
points doubles soient en position élémentaire vis-à-vis de $J$.

\section{Existence de surfaces lagrangiennes dans \vXt}\label{existence}

\subsection{Sections asymptotiquement $J$-holomorphes}\label{sectionsajhol}

Soit \vX\ une variété compacte de dimension $2n$, munie d'un triplet compatible
$(J,\g,\omega)$ avec $\omega$ symplectique et entière ({\it i.e.,} $[\omega]$
provient d'une classe entière en homologie), $E$ un fibré complexe de rang $r$
sur \vX\ muni d'une connexion hermitienne $\nabla^{E}$, $L$ un fibré en droites
holomorphe sur \vX\ tel que $c_{1}(L)=[\omega]\in\HDZ{\vX}$, et $B$ une
connexion hermitienne sur $L$ de courbure $-2i\pi\omega$. (C'est avec ces
signes que $L$ est dit positif par rapport à $J$.) Pour tout $k\in\N^{\ast}$,
$B$ définit une connexion sur $L^{k}$ de courbure $-2ik\pi\omega$, et,
conjointement avec $\nabla^{E}$, elle définit une connexion hermitienne sur
$E\otimes L^{k}$ dont nous noterons $\nabla_k$ la dérivée covariante.

Suivant Donaldson~\cite{Don96} et Sikorav~\cite{Sik98}, nous convenons de
mesurer les dérivées d'une section $s$ de $E\otimes L^{k}$ en munissant \vX\ de
la métrique $\g_{k}=k\g$~; les longueurs mesurées par $\g$ sont multipliées par
$\sqrt{k}$ quand on les mesure avec $\g_{k}$, ceci multiplie la norme des
dérivées $p$-èmes de $s$ et de $\bar{\partial}_{J}s$ par $k^{-p/2}$.

\begin{definitio}
Une suite de sections lisses $s_{k}\in\Gamma(E\otimes L^{k})$ est dite
\textup{asymptotiquement $J$-holomorphe} (en abrégé \textup{$J$-a.h.}) s'il
existe une cons\-tante $C>0$ telle que, pour tout $k\in\N^{\ast}$,
$\|s_{k}\|\leq C$, $\|\bar{\partial}_{J}s_{k}\|_{C^{1}}\leq
\frac{C}{\sqrt{k}}$.
\end{definitio}

\begin{definitio}
Une suite de sections lisses $s_{k}\in\Gamma(E\otimes L^{k})$, $J$-a.h., est
dite \textup{uniformément transverse à la section nulle} si, pour tout $k$
assez grand, il existe une constante $\varepsilon>0$ telle qu'en tout
$\ptx\in\vX$ tel que $\|s_{k}(\ptx)\|\leq\varepsilon$, et que $\|\nabla_k
s_{k}(\ptx)\|\geq\varepsilon$.
\end{definitio}

Donaldson démontre alors le théorème suivant~\cite{Don96}~:

\textit{Soit \vX, $(J,\g,\omega)$, $E$, $L$, $\nabla^{E}$ comme ci-dessus.
Alors il existe une suite de sections $s_{k}\in\Gamma(E\otimes L^{k})$ $J$-a.h.
et uniformément transverse à la section nulle~; pour $k$ assez grand,
$Z_{k}=s_{k}^{-1}(0)$ est une sous-variété symplectique de \vX\ de codimension
$2r$, dont la classe d'homologie est duale à
$(k[\omega])^{r}+\sum_{i=1}^{r}c_{i}(E).(k[\omega])^{r-i}$.}
\newline

Ce résultat est généralisé par Auroux~\cite{Aur97}~:

\textit{Soit \vX, $(J,\g,\omega)$, $E$, $L$, $\nabla^{E}$ comme ci-dessus, et
une suite de sections $\sigma_{k}\in\Gamma(E\otimes L^{k})$ $J$-a.h. Alors,
pour $k$ assez grand, $(\sigma_{k})$ peut être approchée au sens $C^{1}$ par
une suite de sections $s_{k}\in\Gamma(E\otimes L^{k})$ $J$-a.h. et uniformément
transverse à la section nulle.}

\textit{Par ailleurs, si $(J_{t})_{t\in[0,1]}$ est un chemin de structures
presque-complexes sur \vX\ compatibles avec $\omega$, et $(\sigma_{k,t})$ une
suite de chemins de sections $J_{t}$-a.h., dépendant continûment de $t$ ainsi
que leurs dérivées, alors on peut approcher celle-ci au sens $C^{1}$ par une
suite de chemins de sections $s_{k,t}\in\Gamma(E\otimes L^{k})$ $J_{t}$-a.h. et
uniformément transverse à la section nulle~: pour tout $\epsilon>0$, il existe
des constantes $\tilde{K}\geq K$ et $\delta>0$, dépendant de $\epsilon$, de la
géométrie de $\vX$ et des bornes des dérivées des $\sigma_{k,t}$, telles que
pour tout $k>\tilde{K}$, (i) les sections $s_{k,t}$ et leurs dérivées dépendent
continûment de $t$, (ii) pour tout $t\in [0,1]$,
$\|s_{k,t}-\sigma_{k,t}\|<\epsilon$, $\|\nabla_k
s_{k,t}-\nabla_k\sigma_{k,t}\|<\epsilon$, (iii) pour tout $t\in [0,1]$, la
suite $(s_{k,t})$ est $\delta$-uniformément transverse à 0.}
\newline

Signalons aussi le résultat suivant de Donaldson et
Sikorav~\cite{Don96,Sik98}, en termes de courant, dans le cas où $E$
est trivial~:\label{donaldsik}

\textit{Soit $s_{k}\in\Gamma(L^{k})$ une suite de sections $J$-a.h.
et uniformément transverse à la section nulle, et
$Z_{k}=s_{k}^{-1}(0)$ pour tout $k$. Alors la suite des courants
$\frac{1}{k}Z_{k}$ converge vers le dual de Poincaré de $\omega$.
Plus précisément, il existe une constante $C$ telle que pour tout
$\psi\in\Omega^{2n-2}(\vX)$ et tout $k\in\N^{\ast}$,
\[ | \frac{1}{k}\int_{Z_{k}}\psi - \int_{\surX}\omega\wedge\psi |
\leq \frac{C}{\sqrt{k}}\|d\psi\|_{L^{\infty}}, \] la norme
$\|\cdot\|_{L^{\infty}}$ étant prise par rapport à \g.}

\paragraph{}
Dans l'usage que nous allons faire des familles asymptotiquement
holomorphes, il sera pratique d'utiliser aussi les pinceaux de
Lefschetz symplectiques~\cite{Don99} et les revêtements ramifiés
symplectiques de $\CPd$~\cite{Aur00}.

La proposition 6 de~\cite{Don99} implique le résultat suivant~:
\begin{theoreme} \label{existlef}
Pour tous réels $C,\eta>0$ et tout entier $k$ suffisamment grand, il existe des
sections $s_k^{0}$, $s_k^{1}$ de $L^{k}$, $C$-bornées, satisfaisant aux
conditions suivantes~:

$(i)$ $s_k^{0}$ est $\eta$-transverse à la section nulle de $L^{k}$,

$(ii)$ $(s_k^{0},s_k^{1})$ est $\eta$-transverse à la section nulle de
$L^{k}\oplus L^{k}$,

$(iii)$ $\partial(s_k^{1}/s_k^{0})$ est $\eta$-transverse à 0 sur le
complémentaire de $B$, lieu des zéros communs à $s_k^{0}$ et à $s_k^{1}$,

\noindent et telles que $F=s_k^{1}/s_k^{0}$ définisse un faisceau topologique
de Lefschetz de base $B$ pour la structure $k\omega$.
\end{theoreme}

Il est clair, d'après la démonstration de cette proposition par
Donaldson, que pour chaque point $\ptx$ de $\vX$, il existe un tel
pinceau pour lequel $\ptx$ soit un point-base.

Avec les sections $(s^{0}_{k},s^{1}_{k},s^{2}_{k})$ de $\C^{3}\otimes L^{k}$,
D.Auroux parvient à une description des variétés symplectiques comme
revêtements ramifiés du plan projectif complexe~:

\begin{definitio}
Soit $\epsilon>0$, $U$ un voisinage d'un point $\ptx$ dans $\vX$, $\phi~:
U\rightarrow\C^{2}$ une carte locale complexe de classe $C^{1}$, et $J_{0}$ la
structure complexe canonique de $\C^{2}$. L'on dit que $\phi$ est
\textup{$\epsilon$-approximativement holomorphe} si
$|\phi^{\ast}J_0-J|\leq\epsilon$ en tout point de $U$. En d'autres termes, pour
tout vecteur tangent $\vecv$, $|\bar{\partial}\phi(\vecv)| \leq
\frac{\epsilon}{2}|d\phi(\vecv)|$.
\end{definitio}

\begin{definitio}
Une application $f:\vX\rightarrow\CPd$ est dite \textup{localement et
$\epsilon$-holomorphi\-quement modelée en $\ptx$ sur une application
$g:\Cd\rightarrow\Cd$} s'il existe des voisinages $U$ de $\ptx$ dans $\vX$ et
$V$ de $f(x)$ dans \CPd, et des applications $\epsilon$-approximativement
holomorphes $\phi:U\rightarrow\Cd$ et $\psi:V\rightarrow\Cd$, de classe
$C^{1}$, telles que $f_{|U}=\psi^{-1}\circ g\circ\phi$.
\end{definitio}

\begin{definitio}
Une application $f:\vX\rightarrow\CPd$ est un \textup{revêtement singulier
$\epsilon$-approxima\-tivement holomorphe, ramifié au-dessus d'une sous-variété
$R$ de $\vX$}, si $df$ est de rang réel 2 le long de $R$, surjective ailleurs,
et si, en tout $\ptx\in\vX$, elle est $\epsilon$-holomorphiquement modelée sur
l'une des applications suivantes~:

(i) $(z,w)\mapsto (z,w)$,

(ii) $(z,w)\mapsto (z^{2},w)$,

(iii) $(z,w)\mapsto (z^{3}-zw,w)$.

Nous écrirons en abrégé que $f:\vX\rightarrow\CPd$ est un revêtement
$J$-a.h.s. s'il existe un $\epsilon$ tel que $f$ soit un revêtement
singulier $\epsilon$-approximativement holomorphe.
\end{definitio}

(Nous parlerons de surface $\epsilon$-approximativement holomorphe
pour l'image réciproque par un tel revêtement d'une courbe
holomorphe de $\CPd$.)

Rappelons alors le premier théorème obtenu par Auroux dans
\cite{Aur00}~: toute variété symplectique compacte $(\vX,\omega)$ de
dimension 4 peut être réali\-sée symplectiquement comme un
revêtement de \CPd\ ramifié le long d'une courbe symplectique lisse
$R$ dans \vX, se projetant sur \CPd\ en une surface dont les seules
singularités sont des points doubles à croisements normaux et des
points cuspidaux (de première espèce). De plus, pour toute structure
presque-complexe adaptée à $\omega$ sur \vX\ et tout $\epsilon>0$,
nous pouvons demander que l'application de revêtement soit
$\epsilon$-approximativement holomorphe.
\newline

Plus précisément (et en ajoutant un paramètre)~:\label{aurouxparam}

\textit{Soit $(\vX,\omega)$ une variété symplectique compacte de
dimension 4, $(J_{t})_{t\in [0,1]}$ une famille de structures
presque-complexes sur $\vX$, compatible avec $\omega$, $L$ un fibré
en droites holomorphe sur $\vX$ tel que $c_1(L)=[\omega]$, et, pour
tout $k$ entier suffisamment grand, $(\sigma_{k,t})_{k,t\in [0,1]}$
une famille de suites de sections $J_{t}$-a.h. de $\C^{3}\otimes
L^{k}$, telle que les $\sigma_{k,t}$ et leurs dérivées dépendent
continûment de $t$.}

\textit{Alors, pour tout $\epsilon>0$ fixé, il existe (pour $k$
assez grand) une famille $(s_{k,t})_{k\gg 0,t\in [0,1]}$ de suites
de sections $J_{t}$-a.h. de $\C^{3}\otimes L^{k}$, ne s'annulant
nulle part, dépendant continûment de $t$, vérifiant
\[ \forall t\in [0,1], \|s_{k,t}-\sigma_{k,t}\|_{C^{3},\degmetr_{k}}
\leq \epsilon, \] et telle que, pour tout $t$, $s_{k,t}$ définisse
un revêtement ramifié $\vX\rightarrow\CPd$, $J_{t}$-a.h.s.}
\newline

La version à $t$ fixe donne une suite (indexée par l'entier $k$) de
revêtements $\pi_{k}: \vX\rightarrow\CPd$, $J$-a.h.s. (Ceci montre
aussi que l'$\epsilon$ du théorème est en $O(1/\sqrt{k})$.)

D.~Auroux montre comment l'on retrouve des pinceaux à partir de là~:
pour $k$ assez grand, on peut perturber $(s^{0}_{k},s^{1}_{k})$,
extrait du triplet $(s^{0},s^{1},s^{2})$, et obtenir un faisceau.

\subsection{Plans lagrangiens et formes lagrangiennes}\label{plansformeslagr}

L'objet de cette section est d'établir quelques lemmes auxquels nous ferons
appel dans la construction ultérieure.\\

\paragraph{1)}
\begin{lemme} \label{deformeplans}
Soit une collection \textbf{$\Lambda$} de $N$ plans tangents lagrangiens
$\Lambda_{1},\ldots,\Lambda_{N}$ en un point d'un plan complexe hermitien,
transverses deux à deux et en position élémentaire (p.\pageref{poselem}). Il
existe $2N$ plans complexes $\Pi_{1},\Pi_{2},\ldots,\Pi_{2N-1},\Pi_{2N}$ tels
que pour tout $i\in\{1,\ldots,N\}$, $\Pi_{2i-1}$ et $\Pi_{2i}$ soient
transverses à $\Lambda_{i}$ avec des signes d'intersection opposés, et
transverses aux $\Lambda_{j},\quad j\neq i$ avec le même signe d'intersection.
De plus, soit \textbf{$\Pi$} la réunion des $\Pi_{i}$; alors, dans tout
voisinage de \textbf{$\Pi$}, il existe une courbe complexe lisse
désingularisant \textbf{$\Pi$} évitant \textbf{$\Lambda$}.
\end{lemme}
Il suffit de le démontrer dans $\Cd$ muni de sa métrique hermitienne canonique.

Un plan $\Lambda$ d'équation $\{ w=a\bar{z}+bz\}$ dans $\Cd$ est
lagrangien si, et seulement si $|a|^{2}-|b|^{2}=1$. Une droite
complexe d'équation $\{ w=cz\}$ est alors transverse à $\Lambda$
\sii\ $|b-c| \neq |a|$~; le signe de l'intersection est égal à plus
ou moins le signe de $|b-c|-|a|$ selon l'orientation de $\Lambda$.
Deux plans lagrangiens $\Lambda_{1}$ et $\Lambda_{2}$, d'équations
respectives $\{ w=a_{1}\bar{z}+b_{1}z\}$ et $\{
w=a_{2}\bar{z}+b_{2}z\}$, sont transverses \sii\
$\Re(\overline{a_{1}}a_{2}-\overline{b_{1}}b_{2})\neq 1$.

Dans notre cas, $\Lambda_{j}$ a pour équation $w=a_{j}\bar{z}$, avec
$a_{j}=e^{2i(j-1)\pi/N}$. Considérons les plans
$\Pi_{j}=\{(z,w)\in\Cd\tq w=c_{k}z\}$, $k=1,\ldots,2N$, où
$c_{k}=e^{i(k-1)\pi/N}$. La désingularisation de la réunion des
$\Pi_{j}$ retenue est alors la courbe d'équation
\[ w^{2N}-z^{2N}=\varepsilon \hspace{5mm} (\varepsilon>0).
\qed \]

\paragraph{2)}
Nous ferons appel plus tard à un énoncé d'approximation de varié\-tés
lagrangiennes, vues comme des courants, par des formes différentielles
lagrangiennes ({\it i.e.,} orthogonales à $\omegaz$)~:
\begin{lemme} \label{formederham}
Soit $(\vW_{0},\omega_{0})$ une variété symplectique de dimension
$2n$, et $V\subset\vW_{0}$ une sous-variété lagrangienne, compacte,
orientée de $\vW_{0}$. Alors, pour toute métrique riemannienne sur
$\vW_{0}$, pour tout $\varepsilon>0$ et tout voisinage $U$ de $V$
d'ordre $\varepsilon$, il existe une constante $C>0$ et une
$n$-forme différentielle $\alpha$ de classe $\Cinf$ sur $\vW_{0}$,
\textup{fermée}, à support dans $U$, telle que
$\alpha\wedge\omega_{0}=0$, et que, pour toute $n$-forme $\varphi$
sur $\vW_{0}$, on ait
\[ | \int_{\scriptstyle{W}_{0}}\alpha\wedge\varphi - \int_{\surV}\varphi | \leq
C\varepsilon\|\varphi\|_{1}, \] où $\|\cdot\|_{1}$ désigne la norme
$C^{1}$. De plus, le courant que définit $\alpha$ est homologue au
courant $T$ d'intégration sur $V$.
\end{lemme}
{\it Remarque~:} autrement dit, le courant $T$ d'intégration sur $V$
est limite de courants diffus et lagrangiens.
\newline
Suivons le procédé de régularisation du \S 15 du livre de de Rham
\cite{Rha} en l'adaptant au cadre symplectique~:

Premièrement. Sur $\R^{2n}$, pour $\varepsilon<\frac{1}{2}$,
donnons-nous une famille de difféo\-morphismes symplectiques
$s_{\dey}$ ($\vecy\in B^{2n}(0;\varepsilon)$, la boule de centre 0
et de rayon $\varepsilon$ dans $\R^{2n}$), dépendant de manière
\Cinf\ de $\vecy$, qui co\"{\i}ncident avec la famille des
translations par $\vecy$ sur la boule $B^{2n}(0;\varepsilon)$ et
avec l'identité hors de $B^{2n}(0;1)$ (une telle famille se
construit au moyen d'une famille de fonctions génératrices, {\it
cf.}~\cite{McDSa}).

Nous nous donnons aussi une fonction \Cinf\ $f\geq 0$, à support
dans $B^{2n}(0;\varepsilon)$, d'intégrale 1, et nous considérons sur
les formes différentiel\-les l'opération de moyennisation
\[ R^{\ast}\varphi = \int_{\R^{2n}}s_{\dey}^{\ast}(\varphi)\,
f(\vecy)\, d\vecy . \] Cet opérateur agit aussi sur les courants par
la formule
\[ (RT)(\varphi) = T(R^{\ast}\varphi). \]

Les démonstrations des propositions 1 (\cite{Rha}, p.~77) et 2
(\cite{Rha}, p.~78) montrent que $R$ commute à l'opérateur $d$ et
envoie les courants \Cinf\ sur des courants \Cinf\, et les courants
quelconques à support dans $B^{2n}(0;\varepsilon)$ sur des courants
\Cinf.

Deuxièmement. Recouvrons $W_{0}$ au voisinage de $U$ par une collection finie
de cartes de Darboux $\varphi_{i}~: B_{i}(1)\rightarrow W_{0}$ telles que les
$\varphi_{i}(B_{i}(\varepsilon))$ forment encore un recouvrement de $W_{0}$ au
voisinage de $U$. Notons $R_{i}$ l'opérateur de régularisation transporté de
$R$ sur $U_{i}=\varphi_{i}(B_{i}(1))$. \`{A} l'aide d'une partition de l'unité,
nous voyons que le produit $R_{1}\circ R_{2}\circ\ldots\circ R_{N}$ de tous les
opérateurs $R_{i}$, encore noté $R$, satisfait à la propriété suivante~: pour
tout courant $T$ tel que $T\wedge\omega_{0}=0$, $RT$ est de classe \Cinf, à
support dans $U$, et $(RT)\wedge\omega_{0}$ est encore nul.

La forme $\alpha=RT$ ainsi définie est fermée, $R$ commutant avec
$d$. L'estimation sur les normes est celle qu'établit de Rham. \qed

{\it Remarque~:} \label{convergcourant} Du lemme précédent, nous
déduisons que pour toute sous-variété à bord $\tau$, de dimension
$n$, compacte, plongée dans $W$ et transverse à $V$, l'intégrale
$\int_{\tau}\alpha$ tend vers le nombre algébrique d'intersec\-tion
de $V$ avec $\tau$ quand $\varepsilon$ tend vers 0.

Nous aurons aussi besoin du raffinement suivant~:
\begin{lemme} \label{raffine}
En dimension $2n=4$, avec les hypothèses et notations du lemme
précédent, supposons qu'il existe une structure presque-complexe
$J$, compatible avec $\omega$ (telle que $\omega\wedge\omegaz=0$) et
que $V$ soit $J$-holomorphe. Alors on peut s'arranger pour que la
2-forme $\alpha$ donnée par le lemme précédent soit adaptée à $J$,
{\it i.e.,} que $\alpha(J\vecx,\vecy)$ soit une forme symétrique en
$\vecx$ et $\vecy$ le long de $V$.
\end{lemme}
Soit \g\ la métrique compatible déterminée par $J$ et $\omega$, et
$\Omega=\omega\wedge\omega$. Un voisinage (suffisamment petit) de
$V$ dans $\vW_{0}$ s'identifie à un voisinage de la section nulle
dans $T^{\ast}V$ de manière à ce que les fibres soient des droites
$J$-complexes le long de $V$.

Notons $\varsigma$ la symétrie de $T^{\ast}V$ par rapport à $V$ (la
section nulle). Reprenant les notations de la démonstration du
lemme~\ref{formederham}, dans le cas de \Rn, les difféomorphismes
$s_{\dey}$ peuvent être choisis de telle sorte que pour tous
$\vecx$, $\vecy$, nous ayons
$s_{\varsigma\dey}(\varsigma\vecx)=\varsigma(s_{\dey}(\vecx))$, et
la fonction $f$ telle que pour tout $\vecy$,
$f(\varsigma\vecy)=f(\vecy)$. Alors
$R^{\ast}\varsigma^{\ast}=\varsigma^{\ast}R^{\ast}$. Cette propriété
se maintient lorsque $R$ est construit globalement sur $\vW_{0}$ au
voisinage de $U$, comme composée des $R_{i}$.

Soit maintenant $\alpha=RT$, et notons $\vectu$ le champ de
bivecteurs tel que $\alpha=\vectu\prodint\Omega$. La commutation
avec $\varsigma$ assure qu'en chaque point de $V$, $\vectu$ est
invariant par $\varsigma$~; cela signifie que $\vectu$ s'écrit comme
la somme d'un bivecteur tangent à $V$ et d'un bivecteur normal à
$V$. Un tel bivecteur est $J$-invariant par construction, d'où le
lemme. \qed

\paragraph{3)}
Enfin rappelons une notion d'algèbre linéaire utilisée par
Donaldson~\cite{Don96}. Dans $\Cn$ muni de sa métrique euclidienne canonique et
de sa forme symplectique $\omega$ standard, soit $\grass$ la grassmannienne des
$(2n-2)$-plans réels orientés. Sur chaque $\Pi\in\grass$, la métrique
euclidienne canonique définit une forme volume $\Omega_{\Pi}$.
\begin{definitio}
L'\textup{angle de Kähler} $\theta~: \grass\rightarrow [0,\pi]$ est défini pour
tout $\Pi\in\grass$ par
\[ \theta(\Pi)=\arccos
(\frac{1}{(n-1)!}\frac{\omega^{n-1}_{|\Pi}}{\Omega_{\Pi}}). \]
\end{definitio}
Le nombre $\theta$ est une mesure intrinsèque de l'écart aux plans complexes.
L'ensemble des $(2n-2)$-plans orientés $\Pi$ qui sont symplectiques  est égal à
$\theta^{-1}([0,\pi/2[)$.

Ainsi, pour une variété symplectique $(\vX,\omega)$ munie d'une
structure pres\-que-complexe compatible avec $\omega$, et $\vY$ une
sous-variété orientée de codimension 2 de $\vX$, nous pouvons
définir en tout $\pty\in\vY$ un angle $\theta_{\eny}(\vY)$. Pour
revenir à notre situation, Donaldson~\cite{Don96} démontre
l'existence d'une constante $C$ telle que, pour tout $k$, les
$Z_{k}$ donnés par ses théorèmes vérifient
\[ \theta_{\enz}(Z_{k})\leq \frac{C}{\sqrt{k}} \]
pour tout $\ptz\in Z_{k}$.

\subsection{Précisions lagrangiennes}\label{precilagr}

Dorénavant, sauf mention contraire, nous nous placerons en dimension
$2n=4$. Dans le cas particulier où l'on considère des structures
symplectiques orthogonales à $\omegaz$, nous allons chercher à
constru\-ire des surfaces $\omegaz$-lagrangiennes et
$\omega$-symplectiques là où Donaldson et Auroux obtenaient des
surfaces $\omega$-symplectiques, comme lieux d'annulation de
sections asymptotiquement $J$-holomorphes.

\paragraph{} Rappelons le décor~: $(\vXt,\omega)$ est une variété
symplectique fermée de dimension 4, $\omega$ étant à périodes entières~; $L$ un
fibré en droites holomorphe sur \vXt\ tel que $c_{1}(L)=[\omega]\in\HDZ{\vX}$~;
$S$ une surface $\omega$-lagrangienne immergée à croisements normaux~; et
$\omegaz$ une forme symplectique exacte sur $\vXt\setminus S$, orthogonale à
$\omega$ ({\it i.e.,} $\omega\wedge\omegaz=0$).

\paragraph{1)}

\begin{lemme} \label{jholo1}
Donnons nous une suite de sections $\sigma_{k}\in\Gamma(L^{k})$
$J$-a.h. et uniformément transverse à la section nulle, et un réel
$\eta>0$. Il existe une suite, $J$-a.h. et uniformément transverse à
la section nulle, de sections $s_{k}\in\Gamma(L^{k})$,
$C^{0}$-proches de $\sigma_{k}$ à l'ordre $O(1/\sqrt{k})$ et
$\eta$-proches de $\sigma_{k}$ au sens $C^{1}$, telles que, pour $k$
assez grand, $\Sigma_{k}=s_{k}^{-1}(0)$ soit lagrangienne pour
$\omegaz$ sur $\vXt\setminus S$, $J$-holomorphe au voisinage de $S$
et homologue au dual de Poincaré de $k[\omega]$.
\end{lemme}

Commençons par appliquer les résultats de Donaldson et d'Auroux dans le cas où
$E$ est trivial de rang $1$: soit $\sigma_{k}$ une suite de sections  $J$-a.h.
et $\varepsilon$-uniformément transverse à la section nulle, telle que
$Z_{k}=\sigma_{k}^{-1}(0)$ soit uniformément transverse à $S$ pour $k$ assez
grand. (Rappelons que, aussi sur les $Z_{k}$, nous utilisons les métriques
$\g_{k}$.) Pour chaque $k$ assez grand, perturbons $\sigma_{k}$ à l'ordre
$1/\sqrt{k}$ au sens $C^{0}$ et en $o(\varepsilon)$ au sens $C^{1}$ de manière
à rendre $Z_{k}$ transverse à $S$. Toujours pour $k$ assez grand, une
perturbation du même type assure que les plans tangents à $Z_{k}$ aux points de
$Z_{k}\cap S$ sont $J$-complexes, et une seconde perturbation garantit alors
que $Z_{k}$ est $J$-holomorphe (donc $\omegaz$-lagrangienne) au voisinage de
$S$. Notons aussi que la courbure de Gauss et la courbure moyenne de $Z_{k}$
sont bornées indépendamment de $k$~\cite{Sik98}.

Comme $\|\nabla_k\sigma_{k}\|$ est minorée par $\varepsilon$ et que
$\|\nabla_k^{2}\sigma_{k}\|$ est majorée (indépen\-damment de $k$), il existe
un $\delta>0$ tel que $\omegaz$ admet une carte de Darboux sur toute boule de
rayon $\delta$, donnée par une forme de Liouville $\lambdaz$ (primitive de
$\omegaz$), identiquement nulle sur un disque $\omegaz$-lagrangien $D_0$, et
que, pour tout $k$ et tout $\ptz\in Z_{k}$, l'intersection de $Z_{k}$ avec la
boule de centre $\ptz$ et de rayon $\delta$ dans \vXt\ est un disque
correspondant au graphe d'une 1-forme différentielle sur $D_0$.\\

\indent Soit $k > 1/\delta^{2}$, décomposons $Z_{k}$ par des triangles $\tau$
de diamètres inférieurs à $\delta$ (toujours pour la métrique $\g_k$). Et
rendons $J$-holomorphe le plan tangent à $Z_{k}$ aux sommets du complexe formé
par les $\tau$, en ne déformant $Z_k$
qu'à l'ordre $O(1/\sqrt{k})$.\\
\indent Sur chaque triangle $\tau$, nous avons l'égalité
$\int_{\partial\tau}\lambdaz=\int_{\tau}\omegaz$.\\

Démontrons d'abord que l'on peut déformer par isotopie la variété $Z_k$ le long
des arêtes de la triangulation choisie, sans toucher aux plans tangents des
sommets, en ne modifiant les plans tangents le long des arêtes qu'à l'ordre
$O(1/\sqrt{k})$ au sens $C^{0}$ et $\varepsilon$ au sens $C^{1}$, afin
d'annuler les intégrales $\int_{\partial\tau}\lambdaz$ pour tous les
triangles $\tau$.\\

Il existe une 1-forme différentielle $\mu_0$ sur $\vXt\setminus S$ telle que $d
\mu_0=\omegaz$, ce qui assure l'existence d'une constante $C_0$ telle que pour
toute 2-chaîne singulière $T$ sur $Z_k$ on ait
\[
|\int_T\omegaz|=|\int_{\partial T}\mu_0|\leq \frac{C_0}{\sqrt{k}}\|\partial T\|
\]
où $\|\partial T\|$ désigne la longueur de la courbe $\partial T$ pour
$\g_{k}$. Notons $[\omegaz]_\tau$ le 2-cocycle simplicial donné par
l'intégration de $\omegaz$ sur le complexe simplicial $(\tau)$; d'après
~\cite{Sik01}, le théorème de Hahn-Banach fournit une 1-cochaîne simpliciale
$c_0$ pour $(\tau)$ telle que $dc_0=[\omegaz]_\tau$, et telle que pour toute
arête $a$ de $(\tau)$ on ait
\[
|c_0(a)|\leq \frac{C_0}{\sqrt{k}}\|a\|
\]
De plus, comme la courbure de $Z_k$ est bornée, on en déduit
\[
|c_0(a)|\leq \frac{CC_0 \delta}{\sqrt{k}}
\]

En décomposant $\partial\tau=a-b+c$ suivant les trois côtés de $\tau$, on a
\[
\int_{\partial\tau}\lambdaz=c_0(a)-c_0(b)+c_0(c).
\]
Le long de l'arête $a$, considérons une surface $A$ $J_0$-holomorphe tangente à
$a$; puisque $\omegaz$ donne une forme d'aire contrôlée par $g_0$ sur $A$, il
est possible de trouver un chemin $\tilde{a}$ dans $A$ proche de $a$ à l'ordre
$O(1/\sqrt{k})$ au sens $C^{0}$ et $\varepsilon$ au sens $C^{1}$, coïncidant
avec $a$ près des extrémités de sorte que l'aire balayée dans $A$ satisfasse à
\[
\int_{a}\lambdaz-\int_{\tilde{a}}\lambdaz=c_0(a).
\]
On procède de même pour toutes les arêtes. Comme $J$ et $J_0$ sont
orthogonales, les surfaces $A$ et $Z_k$ sont suffisamment transverses pour
qu'on puisse prolonger la déformation des arêtes en une isotopie de $Z_k$, sans
changer la propriété de graphe par rapport à $\lambda_0$. On obtient ainsi de
nouveaux triangles $\tilde{\tau}$.\\
De plus
\[
\int_{\tilde{\tau}}
\omegaz=\int_{\partial\tilde{\tau}}\lambdaz=\int_{\partial\tau}\lambdaz-(c_0(a)-c_0(b)+c_0(c))=0.
\]
\\

Nous déformons ensuite la surface $Z_{k}$ sans modifier davantage le
1-squelette de $(\tilde{\tau})$:

Donnons-nous un triangle $\tilde{\tau}$, et $\pts$ un point de $Z_{k}$
extérieur aux arêtes de $\tilde{\tau}$, centre d'une boule de rayon $\delta$
contenant $\tilde{\tau}$. La surface peut être vue comme le graphe d'une
1-forme $\alpha$ au-dessus d'un ouvert de $\Rd$, et, par nullité de
$\int_{\partial\tilde{\tau}}\lambdaz$, il existe un petit anneau au-dessus
duquel $\alpha$ est exacte, égale à une différentielle $df$. Par ailleurs,
$\alpha$ admet une décomposition sous la forme $d\varphi+d^{\ast}\psi$, unique
si nous imposons que $\varphi=f$ le long des arêtes. Comme, par ailleurs,
$\|d\alpha\|=\|dd^{\ast}\psi\|$ est majorée par $C'/\sqrt{k}$, il s'ensuit que
$\|d\varphi-\alpha\|=\|d^{\ast}\psi\|$ est également majorée par $C"/\sqrt{k}$,
où $C"$ est encore une constante indépendante de $k$.\\
Remplaçons alors le graphe de $\alpha$ à l'intérieur de $\tilde{\tau}$ par
celui de $d\varphi$, ce qui transforme le triangle $\tilde{\tau}$ en un
triangle $\tilde{\tau}'$ lagrangien. Puis modifions encore un peu la surface au
voisinage des arêtes, de manière $\frac{C}{\sqrt{k}}$-proche, afin de la rendre
lisse et $J$-holomorphe sur le 1-squelette de la nouvelle
triangulation.\\

Finalement, il existe une structure presque-complexe $\Jt_k$,
$\frac{C}{\sqrt{k}}$-proche au sens $C^{0}$ de $J$, co\"{\i}ncidant avec $J$
dans un petit voisinage de $S$, et telle que la déformée de $Z_{k}$ reste
$\omegaz$-lagrangienne et $\Jt$-holomorphe. Les estimations sur la déformation
de $Z_{k}$ permettent de déformer les sections $\sigma_{k}$ en des sections
$s_{k}$ possédant les propriétés requises pour le lemme \ref{jholo1}. \qed

\paragraph{2)}

\begin{lemme} \label{jholo2}

Donnons-nous une collection finie de points $\pts_{i}\in S$ et, pour
chaque $i$, une famille finie de droites $J$-complexes $F_{i}^{(j)}$
de $T_{\ens_{i}}W$, deux à deux transverses, ainsi qu'une famille de
suites de sections $\{\sigma_{i,k}^{(j)}\}$ de $L^{k}$, $J$-a.h.,
uniformément transverses à la section nulle, telles que pour tous
$i$ et $j$, $\sigma_{i,k}^{(j)}$ s'annule en $\pts_{i}$, et que
$F_{i}^{(j)}$ soit le plan tangent en $\pts_{i}$ à
$(\sigma_{i,k}^{(j)})^{-1}(0)$. Alors, pour tout $\epsilon>0$, il
existe un entier $K$ tel que pour tout $k\geq K$, il existe une
collection $\{s_{i,k}^{(j)}\}$ de sections de $L^{k}$, $J$-a.h. et
uniformément transverses à la section nulle, $J$-holomorphes au
voisinage de $S$, $\omega$-symplectiques, $\omegaz$-lagrangiennes
sur $\vXt\setminus S$, $C^{0}$-proches des $\sigma_{i,k}^{(j)}$ à
l'ordre $O(1/\sqrt{k})$ et $\epsilon$-proches des
$\sigma_{i,k}^{(j)}$ au sens $C^{1}$ (pour $\g_{k}$), telles que
pour tous $i$ et $j$, $s_{i,k}^{(j)}(\pts_{i})=0$ et $F_{i}^{(j)}$
soit le plan tangent en $\pts_{i}$ à $(s_{i,k}^{(j)})^{-1}(0)$.
\end{lemme}
Cela résulte de l'existence de pinceaux de Lefschetz (théorème~\ref{existlef})
et du lemme~\ref{jholo1}, compte tenu de la possibilité d'imposer les plans
tangents aux points marqués. \qed

\begin{definitio}
L'on dira (de manière laconique) qu'une surface $\Sigma$ est
\textup{$\omega$-holomorphe et $\omegaz$-lagrangienne}, en abrégé
\textup{$(\omega,\omegaz)$-h.l.}, si elle est lagrangienne pour
$\omegaz$ sur $\vXt\setminus S$ et s'il existe une structure
presque-complexe $\Jt$ sur \vXt, compatible avec $\omega$ et
co\"{\i}ncidant avec $J$ au voisinage de $S$, telle que $\Sigma$
soit $\Jt$-holomorphe.
\end{definitio}

\begin{lemme} \label{thelemme}
Soit $\Sigmat$ une surface $(\omega,\omegaz)$-h.l., $\Jt$-holomorphe
dans $\vXt$, qui intersecte $S$ transversalement en $N$ points
réguliers $\pts_{1},\ldots,\pts_{N}$ et en $M$ points singuliers
$\pts_{N+1},\ldots,\pts_{N+M}$. Il existe un réel $r>0$ tel que les
disques $D_{1}$,\ldots, $D_{N}$ de centres les $s_{i}$ ($1\leq i\leq
N$) et de rayon $r$ soient disjoints, et que, pour tout
$(\pts'_{1},\ldots,\pts'_{N})\in\prod_{i=1}^{N}D_{i}$, $\Sigmat$
puisse être déformée en une surface $\Sigmat'$,
$(\omega,\omegaz)$-h.l., vérifiant les propriétés suivantes~:

- $\Sigmat'\cap S$ est exactement égale à
$\{\pts'_{1},\ldots,\pts'_{N},\pts_{N+1},\ldots,\pts_{N+M}\}$,

- les signes d'intersection avec $S$ aux points
$\pts'_{1},\ldots,\pts'_{N}$ sont inchangés,

- les plans tangents aux points $\pts_{N+1},\ldots,\pts_{N+M}$ sont
inchangés.
\end{lemme}
Cela résulte aisément de la paramétrisation des variétés
$\omegaz$-lagran\-giennes $C^{1}$-proches de $\Sigmat$ par des
formes différentielles fermées sur $\Sigmat$. \qed

\section{Disjonction lagrangienne}\label{disjolagr}

Les courbes $\Jt$-holomorphes données par les propositions ci-dessus
ne peuvent, {\it a priori}, être extraites de la partie $W$ de
$\vXt$. Pour constru\-ire à partir de ces courbes une surface
lagrangienne pour $\omegaz$ qui ne rencontre pas $W$, nous allons
adapter une construction de J.~Duval~\cite{Duv}.

\subsection{Lemmes de disjonction locale}\label{disjoloc}

Dans les articles de Lalonde et Sikorav~\cite{LaSi} et de
Polterovich~\cite{Pol} est démontré le résultat suivant~:

\textit{Soient $\Sigma$, $\Sigma'$ deux sous-variétés lagrangiennes, fermées,
orientées d'une variété symplectique \vX\ de dimension 4, qui s'intersectent
trans\-versalement en $m$ points $\ptx_{1},\ldots,\ptx_{m}$, de telle sorte que
l'indice d'inter\-section $\iota_{{\mbox{\scriptsize\sl
x}}_{j}}(\Sigma,\Sigma')$ soit égal à $+1$ pour tout $j$. Alors il existe une
sous-variété (lisse, fermée, orientée) $\tilde{\Sigma}$, lagrangienne,
arbitrairement proche de $\Sigma\cup \Sigma'$, telle que $[\tilde{\Sigma}] =
[\Sigma]+[\Sigma']\in\HdZ{\vX}$.}
\newline

De ce résultat et de~\cite{DoSm}, lemme $(2.10)$ se déduit
facilement le lemme suivant~:
\begin{lemme} \label{deformduval}
Soit $\Sigmat$ une surface immergée $(\omega,\omegaz)$-h.l.
($\Jt$-holomor\-phe) dans \vXt qui intersecte $S$ dans un voisinage $U$ d'un
point régu\-lier \pts\ en exactement un point avec deux plans tangents donnant
des signes d'intersection opposés, et qui coupe $S$ transversalement par
ailleurs. Alors, si $U$ est assez petit, on peut modifier $\Jt$ en une
structure $\Jt'$ compatible avec $\omega$ de manière $C^{0}$-arbitrairement
petite, co\"{\i}ncidant avec $J$ au voisinage de $S$, et déformer par chirurgie
$\Sigmat$ en une surface $\Sigmat'$, $(\omega,\omegaz)$-h.l.
($\Jt'$-holomorphe), ne rencontrant plus $S$ dans $U$, sans changer les
intersections géométriques par ailleurs.
\end{lemme}

{\it Remarque~:} la même construction, jointe au
lemme~\ref{deformeplans}, permet d'éliminer les intersections près
d'un point singulier de $S$ (puisque celui-ci est élémen\-taire au
sens de la définition~\ref{poselem}) lorsqu'en ce point, les plans
tangents de $\Sigmat$ sont convenablement choisis~: le
lemme~\ref{deformeplans} réalise la chirurgie locale de manière
$J$-holomorphe près du point, et l'on recolle grâce au lemme
$(2.10)$ de~\cite{DoSm}.

\subsection{Construction d'une surface lagrangienne disjointe de $S$}\label{lagrdisjos}

Nous allons construire une surface $(\omega,\omegaz)$-h.l., homologue au dual
de Poincaré d'un multiple de $\omega$, et disjointe de $S$.

\'{E}tendant une notion utilisée par J.~Duval dans le cadre des
courbes holomorphes, introduisons d'abord la définition suivante~:

\begin{definitio} Soit $\Sigmat$ une surface $(\omega,\omegaz)$-h.l.
dans \vXt\ qui intersecte $S$ transversalement en $N$ points
réguliers $\pts_{1}$,\ldots,$\pts_{N}$ et en $M$ points singuliers
$\pts_{N+1}$, \ldots,$\pts_{N+M}$~; un nombre réel $r>0$ fixé~; et
$D_{i}$ ($1\leq i\leq N$) le disque de $S$ de centre $\pts_{i}$ et
de rayon $r$ (pour la métrique induite par $\omega$ et $\Jt$). On
dit que $\Sigmat$ est \textup{$r$-souple} si les $D_{i}$ sont
disjoints, à une distance supérieure ou égale à $r$ de tous les
points singuliers de $S$, et si, pour tout
$(\pts'_{1},\ldots,\pts'_{N})\in\prod_{i=1}^{N}D_{i}$, $\Sigmat$
peut être déformée de sorte que $\Sigmat\cap S$ soit exactement égal
à $\{\pts'_{1},\ldots,\pts'_{N},\pts_{N+1},\ldots,\pts_{N+M}\}$, en
gardant $\Sigmat$ $(\omega,\omegaz)$-h.l., et sans modifier les
signes d'intersection dans les $D_{i}$ ($1\leq i\leq N$), ni les
plans tangents aux points $\pts_{i}$ ($N+1\leq i\leq N+M$).
\end{definitio}

Le lemme suivant affirme l'existence de telles surfaces intersectant
$S$ en chaque cellule d'une triangulation assez fine~:
\begin{lemme} \label{rsouplesse}
Pour $m\in\N$ assez grand, il existe un réel $r>0$, une
triangulation $(\tau)$ de $S$ dont les triangles sont de diamètre
inférieur à $r$, et une structure presque-complexe $\Jt$,
$O(1/\sqrt{m})$-proche de $J$, co\"{\i}ncidant avec $J$ au voisinage
de $S$, tels que pour chaque triangle $\tau$, on ait deux surfaces
$(\omega,\omegaz)$-h.l., $2r$-souples, obtenues par des sections
transverses de $L^{m}$, qui coupent ce triangle avec des signes
d'intersection opposés, et, pour chaque point singulier de $S$, un
nombre fini de surfaces $(\omega,\omegaz)$-h.l., $2r$-souples, avec
des droites tangentes $J$-holomorphes imposées en ce point.
\end{lemme}
Partons d'un revêtement ramifié $f_{m}~: \vXt\rightarrow\CPd$ comme
celui que donne le théorème d'Auroux (p.\pageref{aurouxparam}), de
lieu singulier $R_{m}\subset\vXt$ ($R_{m}$ $J$-a.h. plongé). Nous
prenons soin que les singularités de $S$ restent loin de $R_{m}$,
que $R_{m}$ coupe $S$ transversalement, et que
$f_{m}(S\setminus(R_{m}\cap S))$ soit immergé à croisements normaux.

Commen\c{c}ons par considérer un point $\pts_{0}\in R_{m}\cap S$~:
il existe des coordonnées locales $(z_{1},z_{2})$ centrées en
$\pts_{0}$ dans $\vXt$ et des coordonnées locales $(w_{1},w_{2})$
centrées en $f_{m}(\pts_{0})$ dans $\CPd$, telles que
$f_{m}(z_{1},z_{2})=(z_{1}^{2},z_{2})$~; dans ces coordonnées,
$R_{m}$ co\"{\i}ncide avec l'axe $z_{1}=0$. Les coordonnées peuvent
être choisies pour que le plan tangent en $\pts_{0}$ à $S$ ait pour
équation $z_{2}=a\overline{z_{1}}$ avec $|a|=1$.

Les deux coniques d'équations $w_{1}=4w_{2}^{2}$ et
$w_{1}=\frac{1}{4}w_{2}^{2}$ se relèvent par $f_{m}$ en deux
surfaces symplectiques $\Sigma'$ et $\Sigma''$, qui donnent près de
$\pts_{0}$ quatre branches $z_{1}=\pm 2z_{2}$ et $z_{2}=\pm 2z_{1}$,
dont deux coupent $S$ positivement et deux négativement.

Une petite perturbation (lemme~\ref{thelemme}) et une
désingularisation hors de $S$ fournissent une surface $\Sigma_{0}$
qui coupe $S$ près de $\pts_{0}$ en deux points $\pts'_{0}$ et
$\pts''_{0}$, avec, en chacun de ces points, des branches
d'intersections positive et négative. Il existe $r_{0}>0$ tel que
$\Sigma_{0}$ soit $r_{0}$-souple~; pour chaque $\pts_{0}\in
R_{m}\cap S$, nous considérons deux triangles équilatéraux sur $S$,
de diamètre $r_{0}$, centrés respectivement en $\pts'_{0}$ et
$\pts''_{0}$, ayant un côté commun centré en $\pts_{0}$.

Dans la suite, $U$ désigne la réunion des losanges sur $S$ centrés
aux points de $R_{m}\cap S$, de diamètre $2r_{0}$, construits comme
ci-dessus.

\vspace{2mm}

Pour tout $\theta\in\CPdual=(\CPd)^{\ast}$, soit $D_{\theta}$ la droite
complexe de $\CPd$ définie par $\theta$, et $\Sigmamth$ l'image réciproque de
$D_{\theta}$ par $f_{m}$. Chaque surface $\Sigmamth$ est symplectique et lisse
hors de $\Sigmamth\cap R_{m}$.

Notons $\Transv^{(m)}$ (respectivement $\Hs^{(m)}$) l'ensemble des
$\theta\in\CPdual$ pour lesquels $\Sigmamth$ n'est pas transverse à
$S$ (respectivement, passe par un point $\pts$ donné dans $S$). Pour
chaque point $\pts$ de $S\setminus U$, considérons l'application
$\tau_{\ens}$ de $S\setminus\{\pts\}$ dans $\CPdual$ qui à $\ptt$
associe la droite $(f_{m}(\pts)f_{m}(\ptt))$ lorsque
$f_{m}(\pts)\neq f_{m}(\ptt)$, et la tangente à la branche de
$f_{m}(\ptt)$ dans $f_{m}(S)$ lorsque $f_{m}(\pts)=f_{m}(\ptt)$.
L'intersection de $\Transv^{(m)}$ avec $\Hs^{(m)}$ est l'ensemble
des valeurs critiques de $\tau_{\ens}$.

Par le théorème de Sard, $\Transv^{(m)}\cap\Hs^{(m)}$ est de mesure
nulle dans $\Hs^{(m)}$ ~; $\Transv^{(m)}$ possède donc une base de
voisinages dans $\CPdual$ dont l'intersection avec $\Hs^{(m)}$ peut
être rendue de mesure aussi petite que l'on veut, uniformément en
$\pts$ sur $S\setminus U$.

Définissons, en chaque $\pts\in S$, $\Hs^{+}$ (respectivement
$\Hs^{-}$) comme le sous-ensemble des points de $\Hs^{(m)}$ où
$\Sigmamth$ coupe $S$ positivement (resp.~négative\-ment). Comme $S$
est totalement réelle, ces deux espaces sont de mesure non nulle,
minorée uniformément sur $S\setminus U$~; nous pouvons donc trouver
un voisinage $\Vtransv$ de $\Transv^{(m)}$ dans $\CPdual$, de mesure
assez petite pour que $\Transv^{(m)}\setminus\Vtransv$ rencontre
$\Hs^{+}$ et $\Hs^{-}$ en tout $\pts\in S\setminus U$. En dehors de
$\Vtransv$, nous disposons de bornes {\it a priori} sur le nombre de
points de $\Sigmamth\cap S$ et sur le minimum de leurs distances
mutuelles. En vertu des lemmes~\ref{jholo1} et \ref{thelemme}, il
existe un réel $r>0$ tel que $\Sigmamth$ soit $r$-souple sur
$S\setminus U$.

Pour conclure, choisissons une triangulation de maille inférieure à
$r$ et à $r_{0}$ pour $U$, ainsi que des surfaces $\Sigmamth$ et
$\Sigma_{0}$ pour $U$ en nombre suffisant. Il suffit alors de
déformer le tout avec le lemme~\ref{jholo2} pour obtenir des
surfaces $(\omega,\omegaz)$-h.l. satisfaisant aux propriétés du
lemme. \qed

\noindent De plus, nous pouvons supposer qu'aux points singuliers de $S$, les
surfaces données par le théorème soient tangentes à des droites $J$-complexes
désingularisables au sens du lemme~\ref{deformeplans}. Appelons $\Sigma_{1}$ la
réunion de toutes ces surfaces.

\vspace{2mm}

Pour chaque triangle $\tau$ et pour $\Sigma_{j}$ ($j\in I$) une surface fermée
orientée immergée dans $\vXt$ intersectant $S$ transversalement, nous noterons
$p_{j}(\tau)$ (respectivement $n_{j}(\tau)$) le nombre de points d'intersection
positive (respectivement négative) de $\Sigma_{j}$ avec $\tau$, et
$i_{j}(\tau)=p_{j}(\tau)-n_{j}(\tau)$. C'est un 2-cocycle car $S$ est fermée de
dimension $2$, et c'est un 2-cobord dès que les classes d'homologie de
$\Sigma_{j}$ et de $S$ ont un nombre total d'intersection égal à $0$; ce qui
sera notre cas puisque l'intégrale de $\omega$ sur $S$ est nulle et que
$[\Sigma_{j}]$ sera un multiple entier de $\omega$. Par ailleurs, pour
$x\in\R$, nous noterons $x^{+}=\sup(x,0)$ et $x^{-}=\sup(0,-x)$.

Une première déformation de $\Sigma_{1}$ permet (par $r$-souplesse)
d'élimi\-ner suffisamment de paires de points d'intersection avec
$S$ de signes opposés pour se ramener à une surface
$(\omega,\omegaz)$-h.l. $\Sigma_{2}$ pour laquelle $p_{2}=i_{1}^{+}$
et $n_{2}=i_{1}^{-}$, en assurant (par le lemme~\ref{deformduval})
qu'il existe une structure presque-complexe $J_{2}$ orthogonale à
$\omegaz$, compatible avec $\omega$, $C^{0}$-arbitrai\-rement proche
du $J$ initial, et qui co\"{\i}ncide avec $J$ au voisinage de $S$,
pour laquelle $\Sigma_{1}$ et $\Sigma_{2}$ sont $J_{2}$-holomorphes.

Fixons $\epsilon>0$. D'après le lemme~\ref{formederham}, il existe
une 2-forme fermée $\omega_{2}$, orthogonale à $\omegaz$, qui
approche la surface $\Sigma_{2}$ à $\epsilon$ près. Il existe alors
un entier $A$ assez grand pour que la forme
$\omega_{3}=A\omega-\omega_{2}$, orthogonale à $\omegaz$, soit elle
aussi symplectique et induise une forme de signe constant sur chaque
triangle $\tau$, d'intégrale $-i_{1}(\tau)$ ({\it cf.} la remarque
p.~\pageref{convergcourant}). Comme $\omega$ et $\omega_{2}$,
$\omega_{3}$ est à périodes entières.

D'après le lemme~\ref{raffine}, il existe une structure
presque-complexe $J_{3}$, compatible avec $\omega_{3}$, orthogonale
à $\omegaz$, $C^{0}$-arbitrairement proche du $J$ initial, telle que
$J_{2}$ et $J_{3}$ co\"{\i}ncident le long de $\Sigma_{1}$ et
$\Sigma_{2}$.

\vspace{2mm}

Nous avons maintenant besoin d'un résultat d'approximation des
formes symplectiques par les courants, qui opère une
``moyennisation'' de la proposition de Donaldson et
Sikorav~\pageref{donaldsik}. Il remplacera la formule de Crofton
dans la preuve de J.~Duval~\cite{Duv}. Dans ce paragraphe, $\omega$
désigne une forme symplectique à périodes entières quelconques sur
la variété $\vX$.

Pour $\varepsilon>0$ ($\varepsilon\geq C/\sqrt{k}$, où $C$ est une constante
indépendante de $k$), donnons-nous comme dans~\cite{Aur00} un triplet de
sections $(s_{k}^{0},s_{k}^{1},s_{k}^{2})$ de $L^{k}$ qui définit un revêtement
ramifié singulier $\varepsilon$-approximativement holomorphe
$f_{k}:\vX\rightarrow\CPd$, avec une surface de ramification $R_{k}\subset\vX$.

Pour chaque point $\theta$ du plan dual $\CPdual=(\CPd)^{\ast}$, considérons la
droite complexe $D_{\theta}\subset\CPd$ définie par $\theta$ et
$f_{k}^{-1}(D_{\theta})=\Sigmakth\subset\vX$. Chaque surface $\Sigmakth$ est
symplectique, lisse en-dehors de $R_{k}\cap\Sigmakth$,
$\varepsilon$-approximativement holomorphe~; elle est le lieu des zéros d'une
section $\skth$ de $L^{k}$. L'holomorphie approximative de $f_{k}$, ainsi que
sa structure le long de $R_{k}$, entra\^{\i}nent que dans chaque pinceau
$\{\theta\in(\CPu)^{\ast}\subset\CPdual\}$, seul un nombre fini de $\Sigmakth$
possède des points singuliers.

Considérons le courant $D_{k}$ de degré 2 sur $\vX$ obtenu en
intégrant les courants $[\Sigmakth]$ pour la mesure de Fubini-Study
$d\theta$ sur $\CPdual$, normalisée de manière à ce que la masse
totale de $\CPdual$ soit égale à 1.

Pour l'analyse de $D_{k}$, nous utilisons encore les métriques $\g_{k}=k\g$ sur
$\vX$. Comme pour le lemme~\ref{jholo1}, il existe un réel $\delta>0$ tel que
les intersections, hors de $R_{k}$, des $\Sigmakth$ avec les boules de rayon
$\delta$ soient des graphes sur leur plan tangent avec les dérivées premières
bornées par 1. La restriction de $D_{k}$ à chaque boule $B(\ptx_{0},\delta)$
qui évite $R_{k}$ est un courant lisse, {\it i.e.,} induit une forme
différentielle $\omega_{(k)}$ sur cette boule.

\begin{lemme} \label{approcheforme}
$\| \omega_{(k)}-k.\omega \| = O(\sqrt{k})$.
\end{lemme}
Suivant Donaldson (\cite{Don96}, p.~702), si $\skth$ est la section
de $L^{k}$ définissant $\Sigmakth$, la 1-forme de connexion associée
$\Akth=(\skth)^{-1}\nabla\skth$ est $L^{1}$, et l'on a
$d\Akth=\Sigmakth-k\omega$. Il existe $\eta>0$ tel que $\skth$ est
$\eta$-transverse à la section nulle. Pour établir le lemme, il
suffit donc d'estimer le courant lisse
\[ B_{k}=\int_{H} d\theta\,d\Akth .\]

Prenons un 2-simplexe singulier $\Delta$ dans $B(\ptx_{0},\delta)$ et
découpons-le en triangles $\Delta_{j}$ de diamètres respectifs $\eta_{j}$.
Notons $\mathbf{C}_{j}$ l'ensemble des $\theta\in\CPdual$ tels que $\|
\skth(\ptx) \| \geq \eta_{j}$ pour tout $\ptx\in\Delta_{j}$, et
$\mathbf{C}_{j}^{\ast}$ son complémentaire dans $\CPdual$. Alors~:
\[ \int_{\Delta_{j}} \!\!\! B_{k} =
\int_{\partial\Delta_{j}}\!\!\! dz\:\int_{\mathbf{C}_{j}}\!\!\!
d\theta\,\frac{\nabla\skth}{\skth} + \int_{\partial\Delta_{j}}\!\!\!
dz\:\int_{\mathbf{C}_{j}^{\ast}}\!\!\! d\theta\,\frac{\nabla\skth}{\skth}.
\]

La première intégrale est nulle, comme somme d'une forme fermée.\\
La seconde est majorée par $C_{1}.\|\partial\Delta_{j}\|.\eta_{j}.\sqrt{k}$ (où
$\|\partial\Delta_{j}\|$ est la lon\-gueur du périmètre du triangle
$\Delta_{j}$, et $C_{1}$ une constante positive indépen\-dante de $k$). En
effet, il existe une constante $C_{2}>0$ (indépendante de $k$) telle que
$\|\nabla\skth\| \leq C_{2}$.\\

Quant à $\int_{C_{j}^{\ast}}d\theta\,/\,\|\skth\|$, montrons qu'elle est
majorée par $C_{7}.\eta_{j}.\sqrt{k}$, où $C_{7}$ est une constante positive
indépendante de $k$~: \\

Par définition de $\mathbf{C}_{j}$, pour tout $\theta\in
\mathbf{C}_{j}^{\ast}$, il existe un point $\ptx\in\Delta_{j}$ tel que
$\|\skth(\ptx)\|<\eta_{j}$. Comme, par ailleurs, les dérivées premières des
$\skth$ sont bornées par 1 sur $B(\ptx_{0},\delta)$, et que le diamètre de
$\Delta_{j}$ est $\eta_{j}$, nous obtenons l'existence d'une constante $C_3>0$
(indépendante de $j$ et de $k$) telle que
\[ \forall\theta\in C_{j}^{\ast},
\forall\ptz\in\partial\Delta_{j}, \|\skth(\ptz)\|<\kappa.\eta_{j}.
\] Quitte à réduire les $\eta_{j}$, nous pouvons supposer que
$\kappa.\eta_{j}<\eta$.

Soit $\ptz\in\partial\Delta_{j}$~; dans $\CPd$, il existe des coordonnées
locales $(z,w)$ centrées en $f_{k}(\ptz)$, telles que la mesure $d\theta$ se
décompose en $d\beta\,du$, où $\{w=\beta z+u\}$ est l'équation locale de
$D_{\theta}$ (en particulier, $\theta\in H_{\enz}^{(k)}$ \sii\ $u=0$). Par
$\eta$-transversalité de $\skth$, et en ayant préalablement réduit les
$\eta_{j}$ au besoin pour rester dans les ouverts de coordonnées locales, il
existe une constante $\kappa'>0$ telle que, pour tout $\theta=(\beta,u)\in
\mathbf{C}_{j}^{\ast}$, $u$ appartienne au disque
$d_{\kappa'\eta_{j}}\subset\C$ de centre 0 et de rayon $\kappa'.\eta_{j}$.

Nous obtenons alors une majoration sur l'intégrale~:
\[ \int_{\mathbf{C}_{j}^{\ast}}\frac{d\theta}{\|\skth\|} \leq
\int_{H_{\enz}^{(k)}}\!\!\!d\beta\:\int_{d_{\kappa'\eta_{j}}}\frac{du}{\|s_{k}^{(\beta,u)}(z)\|}.
\]
Mais $\|s_{k}^{(\beta,u)}(z)\|$ est minoré par $C_4u/\sqrt{k}$ lorsque $u$ est
assez petit, par $\eta$-transversalité ($C'$ cons\-tante), de sorte qu'il
existe une constante $C_5$ telle que
\[ \int_{d_{\kappa'\eta_{j}}}\frac{du}{\|s_{k}^{(\beta,u)}(z)\|} \leq
C_5.\eta_{j}.\sqrt{k}.\] Par ailleurs, $\int_{H_{\enz}^{(k)}}d\beta$ est
uniformément bornée par une constante $C_6$, et il suffit de poser
$C_{7}=C_5C_6$ pour obtenir la majoration souhaitée.

Pour tout triangle $\Delta_{j}$, $|\partial\Delta_{j}|$ est en
$O(\eta_{j})$, de sorte que $\int_{\Delta_{j}} \!\!\! B_{k}$ est
majoré par $C_{4}.\eta_{j}^{2}.\sqrt{k}$ ($C_{4}$ constante
indépendante de $j$ et de $k$). En sommant sur tous les
$\Delta_{j}$, nous obtenons finalement l'existence d'une constante
$C_{5}$ telle que $|\int_{\Delta}B_{k}| \leq C_{5}\sqrt{k}$, ce qui
démontre le lemme. \qed

\vspace{2mm}

Cette construction nous sert à établir le lemme suivant~:
\begin{lemme} \label{duv3}
Soit $K$ un nombre réel strictement positif. Pour tout $\epsilon>0$
suffisamment petit, et pour tout entier $k$ suffisamment grand, il existe un
entier $l=l(k)$, un fibré $L'$ sur $\vXt$, et une famille de sections $s_{3}$
de $L^{k}\otimes L'$, asymptotiquement $J_{3}$-holomorphe et
$\epsilon$-transverse à $\Sigma_{2}$, tels que $\Sigma_{3}=s_{3}^{-1}(0)$ soit
$(\omega,\omegaz)$-h.l. et vérifie
\[ \| n_{3}-l.i_{1}^{+} \| + \| p_{3}-l.i_{1}^{-} \| < \frac{l}{4K} \]
\end{lemme}

C'est maintenant à la structure $\omega_{3}$ que nous allons appliquer la
con\-struction précédente, avec $k>C^{-2}/\epsilon^{2}$~; cette fin de
démonstration est analogue à celle du lemme 3 de~\cite{Duv}.

Donnons-nous comme précédemment $f_{k}:\vXt\rightarrow\CPd$,
$\epsilon$-approximative\-ment holomorphe, et soit
$R_{k}\subset\vXt$ le lieu de ramification (asymptotiquement
$J_{3}$-holomorphe). Comme dans la démonstration du lem\-me
\ref{rsouplesse}, commen\c{c}ons par construire de petits losanges
$U'$ autour des points de $R_{k}\cap S$, qui évitent $\Sigma_{2}$.
Notons à nouveau $U$ la réunion de ces losanges, $\Transv^{(k)}$
(respectivement $\Hs^{(k)}$) l'ensemble des $\theta\in\CPdual$ pour
lesquels $\Sigmakth$ n'est pas transverse à $S$ (respectivement,
passe par un point $\pts$ donné dans $S$). Il existe un voisinage
$\Vtransv$ de $\Transv^{(k)}$ tel que la mesure de
$\Vtransv\cap\Hs^{(k)}$ soit majorée par $\epsilon$ uniformément
pour $\pts\in S\setminus U$, et un réel $r'>0$ tel que pour tout
$\theta\in\Hs^{(k)}\setminus\Vtransv$, la surface $\Sigmakth$ soit
$r'$-souple. Choisissons une triangulation $(\tau')$ plus fine que
$(\tau)$ en triangles de
diamètre inférieur ou égal à $r'$.\\
De plus, soit $\Transv'\,^{(k)}$ l'ensemble des $\theta\in\CPdual$
tels que $\Sigmakth$ ne coupe pas $\Sigma_{2}$ transversalement {\it
et} positivement. La surface $\Sigma_{2}$ étant $J_{2}$-holomorphe
avec $J_{2}$ $\epsilon$-proche de $J$, $f_{k}$ est encore
approximativement holomorphe à l'ordre $O(\epsilon)$ pour $J_{2}$ au
voisinage de $\Sigma_{2}$~; comme chaque $\Sigmakth$ est image
réciproque par $f_{k}$ de la droite holomorphe $D_{\theta}$ de
$\CPd$, la mesure de $\Transv'\,^{(k)}$ est majorée par $\epsilon$
fois une constante $C'$.

Notons $\Vtransv'$ un voisinage de $\Transv'\,^{(k)}$ de mesure inférieure ou
égale à $2C'\epsilon$, et $\Vtransv''$ la réunion de $\Vtransv$ et de
$\Vtransv'$. Alors, en utilisant la remarque p.~\pageref{convergcourant}, pour
chaque triangle $\tau'$, la différence $\int_{\tau'}\omega_{3} -
\int_{H\setminus V''}(\Sigmakth\cdot\tau')\,d\theta$ est (au plus) de l'ordre
de $C''.\epsilon.\mbox{aire}(\tau')$ (où $C''$ est une constante indépendante
de $k$).

Considérons maintenant un entier $n$ et $n$ élements $\theta_{j}$ de $\CPdual$
répartis régulièrement par rapport à la mesure $d\theta$ dans $\CPdual$, et
notons $\Sigma_{(n)}$ la réunion des surfaces $\Sigma_{k}^{(\theta_{j})}$ pour
tous les $\theta_{j}\notin\Vtransv''$. Par convergence des sommes de Riemann
vers l'intégrale et d'après le lemme \ref{approcheforme}, $n$ peut être choisi
assez grand pour que les différences
$\frac{1}{kn}i_{(n)}-\int_{\tau'}\omega_{3}$ soient de l'ordre de
$\epsilon.\mbox{aire}(\tau')$ pour tout triangle $\tau'$. Par $r'$-souplesse et
d'après le lemme~\ref{thelemme}, nous pouvons déformer $\Sigma_{(n)}$ en une
courbe $J_{3}$-holomorphe $\Sigma_{3}$ qui vérifie $p_{3}=i_{(n)}^{+}$ et
$n_{3}=i_{(n)}^{-}$ sur $(\tau')$. La transversalité à $\Sigma_{2}$ a déjà été
assurée.

Finalement, $\omega_{3}$ étant $\epsilon$-proche d'une forme de signe constant
sur chaque $\tau$, il existe une constante $C$ ne dépendant que de la
triangulation telle que $\| n_{3}-l.i_{1}^{+} \| + \| p_{3}-l.i_{1}^{-} \| <
C.l.\epsilon$, avec $l=kn$. Si $\epsilon$ a été choisi suffisamment petit au
départ, le lemme en résulte. \qed\\

Rappelons que la triangulation $(\tau)$ de $S$ est fixée; par linéarité de
l'application cobord de l'espace des 1-cochaînes simpliciales à coefficients
réels $\mathcal{C}^{1}(S)$ dans l'espace des 2-cocycles simpliciaux à
coefficients réels $\mathcal{Z}^{2}(S)$, il existe une constante strictement
positive $K$, telle que, pour tout 2-cocycle $u$, nul en cohomologie, on puisse
trouver une 1-cochaîne $c$ telle que $d(c)=u$ et $\|c\|\leq K\|u\|$. C'est à
cette constante $K$ que nous allons appliquer le lemme précédent.\\

\begin{lemme} \label{duv4}
Pour $k$ entier assez grand, il existe une famille de sections $s_{4}$ de
$L^{k}$ tels que $\Sigma_{4}=s_{4}^{-1}(0)$ soit $(\omega,\omegaz)$-h.l.,
$r$-souple, à intersections positives avec $\Sigma_{3}$, et vérifie
$n_{4}(\tau)>p_{3}(\tau)$, $p_{4}(\tau)>n_{3}(\tau)$, et
$i_{4}(\tau)=-i_{3}(\tau)$ pour tout triangle $\tau$.
\end{lemme}
Partons cette fois de la réunion $\Sigma_{4}'$ de $l$ surfaces
lagrangiennes vérifiant les mêmes propriétés que celles que nous
avions demandées pour $\Sigma_{1}$ ci-dessus. Par le
lemme~\ref{duv3}, $li_{1}+i_{3}$ est le cobord d'une 1-cocha\^{\i}ne
$c$ vérifiant $\|c\|<\frac{l}{4}$. Fixons une arête orientée $a$ de
la triangulation $(\tau)$~; $S$ étant orientée près de $a$, l'un des
deux triangles adjacents à $a$ induit la bonne orientation sur $a$.
Il est alors possible de faire passer $c(a)$ intersections orientées
de ce triangle vers l'autre~: en effet, pour chaque triangle $\tau$,
$\Sigma_{4}'$ intersecte positivement (resp. négativement) $\tau$ en
au moins $l$ points, alors que
$\sum_{\partial\tau}|c(a)|<\frac{3l}{4}$. La courbe $\Sigma_{4}$
ainsi obtenue convient. \qed

Nous pouvons maintenant déformer $\Sigma_{4}$ dans chaque $\tau$ de manière à
faire co\"{\i}ncider chaque point de $\Sigma_{3}\cap\tau$ avec un point de
$\Sigma_{4}\cap\tau$ de signe d'intersection opposé. Par déformation, nous
pouvons aussi faire co\"{\i}ncider les points restants de $\Sigma_{4}\cap\tau$
par paires de points pour lesquels les signes d'intersection sont opposés.
Appliquant le lemme~\ref{deformduval} à la réunion de $\Sigma_{3}$ et de la
déformée de $\Sigma_{4}$, nous éliminons toutes ses intersections avec $S$. Les
intersections de $\Sigma_{3}$ et $\Sigma_{4}$ sont toutes transverses et
positives~: elles peuvent donc aussi être éliminées par chirurgie, d'où le
lemme suivant~:
\begin{lemme}
Il existe une structure symplectique $\omega'$ cohomologue à
$\omega$, et une surface $\Sigmat$ $\omega'$-holomorphe et
$\omegaz$-lagrangienne, disjointe de $S$, homologue au dual de
Poincaré d'un multiple de $\omega$.
\end{lemme}

\section{Dernière étape~: disjonction de $W$}\label{disjow}

D'après Andreotti et Frankel~\cite{AnFr,GrHa}, le champ de Liouville
$\vectxi_{0}$ de $\omegaz$ sur $W$ se prolonge à $\vXz$ en un champ
de Liouville. En suivant les lignes de gradient de $-\vectxi_{0}$
(qui est complet), nous pouvons alors tirer $\Sigmat$ en une surface
$\Sigma$ plongée dans $\vXz$, homologue au dual de Poincaré d'un
multiple de $\omega$. Ceci termine la démonstration du
théorème~\ref{theotrois}.

\newpage

\end{document}